\documentclass[11pt]{article}

\usepackage{amsmath}
\usepackage{amsthm}
\usepackage{amssymb}
\usepackage{amscd}

\theoremstyle{definition}
\newtheorem{theorem}{Theorem}[section]
\newtheorem{prop}[theorem]{Proposition}
\newtheorem{lemma}[theorem]{Lemma}
\newtheorem{corollary}[theorem]{Corollary}
\newtheorem{definition}[theorem]{Definition}

\newtheorem{remark}[theorem]{Remark}
\newtheorem{conjecture}[theorem]{Conjecture}
\newtheorem{question}[theorem]{Question}
\newtheorem{problem}[theorem]{Problem}

\numberwithin{equation}{section}
\newenvironment{demo}[1]{%
  \trivlist
  \item[\hskip\labelsep
        {\it #1.}]
}{%
\hfill\qedsymbol
  \endtrivlist
}

\newcommand\Nat{\mathbb{N}}
\newcommand\Int{\mathbb{Z}}
\newcommand\Rat{\mathbb{Q}}
\newcommand\Comp{\mathbb{C}}

\newcommand\vectorone{\boldsymbol{1}}
\newcommand\PF{\operatorname{PF}}

\newcommand\Cat{\operatorname{Cat}}

\newcommand\type{\operatorname{type}}

\newcommand\GL{\mathbf{GL}} 
\newcommand\Sym{{\mathfrak{S}}} 
\newcommand\Par{\mathcal{P}}
\newcommand\ch{\operatorname{ch}}
\newcommand\sh{\operatorname{sh}}
\newcommand\Res{\operatorname{Res}}
\newcommand\Ind{\operatorname{Ind}}
\newcommand\Hom{\operatorname{Hom}}

\newcommand\Irr{\operatorname{Irr}}
\newcommand\gauss[3]%
     {\bigg[\genfrac{}{}{0pt}{0}{#1}{#2}\bigg]_{#3}}

\newcommand\triv{\operatorname{triv}}

\newcommand{\vectx}{\boldsymbol{x}}
\newcommand{\vecty}{\boldsymbol{y}}

\renewcommand\tilde{\widetilde}
\renewcommand\hat{\widehat}
\newcommand\ep{\varepsilon}

\allowdisplaybreaks

\title{
On the Existence of Generalized Parking Spaces
\\
for Complex Reflection Groups
}

\author{
Yosuke ITO\footnote{
current address : Meiji Yasuda Life Insurance Company,
1-1, Marunouchi 2-chome, Chiyoda-ku, Tokyo 100-0005, Japan
}
, \ 
and 
Soichi OKADA
\\
Graduate School of Mathematics, Nagoya University
\\
Furo-cho, Chikusa-ku, Nagoya 464-8602, Japan
}

\date{
}

\begin{document}

\maketitle

\begin{abstract}
Let $W$ be an irreducible finite complex reflection group acting 
on a complex vector space $V$. For a positive integer $k$, we 
consider a class function $\varphi_k$ given by $\varphi_k(w) = k^{\dim V^w}$ 
for $w \in W$, where $V^w$ is the fixed-point subspace of $w$. 
If $W$ is the symmetric group of $n$ letters and $k=n+1$, then 
$\varphi_{n+1}$ is the permutation character on (classical) 
parking functions. In this paper, we give a complete answer 
to the question when $\varphi_k$ (resp. its $q$-analogue) is 
the character of a representation (resp. the graded character 
of a graded representation) of $W$. As a key to the proof 
in the symmetric group case, we find the greatest common divisors 
of specialized Schur functions. And we propose a unimodality 
conjecture of the coefficients of certain quotients of principally 
specialized Schur functions.
\end{abstract}

\section{%
Introduction
}

\subsection{%
Background
}

A (classical) \emph{parking function} of length $n$ is a map 
$f : \{ 1, 2, \dots, n \} \to \{ 1, 2, \dots, n \}$ satisfying 
$\# f^{-1}( \{ 1, 2, \dots, k \} ) \ge k$ for each $k = 1, \dots, n$.
The notion of parking function was introduced by Pyke \cite{Pyke} and 
independently by Konheim--Weiss \cite{KW}, and 
later becomes one of the main characters in algebraic combinatorics.

Let $\PF_n$ be the set of all parking functions of length $n$.
The symmetric group $\Sym_n$ of $n$ letters $\{ 1, 2, \dots, n \}$ acts 
on $\PF_n$ by the rule $f \mapsto f \circ w$ for $f \in \PF_n$ 
and $w \in \Sym_n$.
The corresponding permutation module is called the (classical) \emph{parking space}.
It is an important problem in Catalan combinatorics to generalize the parking space 
from $\Sym_n$ to other complex reflection groups 
and from the Coxeter number to more general parameters (including the Fuss cases).

It is known (\cite[Proposition~2.6.1]{Haiman}) that the $\Sym_n$-action on $\PF_n$ 
is isomorphic to the $\Sym_n$-action on $(\Int_{n+1})^n / \langle \vectorone \rangle$, 
where $\Int_{n+1} = \Int/(n+1)\Int$ and $\langle \vectorone \rangle$ is the subgroup 
generated by $\vectorone = (1, 1, \dots, 1) \in (\Int_{n+1})^n$.
Hence we see that the corresponding permutation character $\chi_{\PF_n}$ is given by
\begin{equation}
\label{eq:char_PF}
\chi_{\PF_n}(w) = (n+1)^{l(\type(w))-1}
\quad(w \in \Sym_n),
\end{equation}
where $\type(w)$ is the cycle type of $w$ 
and $l(\type(w))$ is the number of cycles in $w$.
In particular, there are $(n+1)^{n-1}$ parking functions and 
the number of $\Sym_n$-orbits is equal to the $n$th Catalan number.

As a generalization of the $\Sym_n$-set $\PF_n$ to Weyl groups $W$, 
Haiman \cite{Haiman} studied the quotient $Q/(h+1)Q$ of the root lattice $Q$, 
where $h$ is the Coxeter number of $W$.
For real reflection groups $W$, Armstrong--Reiner--Rhoades \cite{ARR} 
proposed two new $W$-parking spaces $\Comp[{\mathsf{Park}}^{\text{NC}}_W]$ and 
${\mathsf{Park}}^{\text{alg}}_W$. 
And Rhoades \cite{Rhoades} gives Fuss analogues $\Comp[{\mathsf{Park}}^{\text{NC}}_W(p)]$ and 
${\mathsf{Park}}^{\text{alg}}_W(p)$ of the constructions in \cite{ARR}.
Note that ${\mathsf{Park}}^{\text{NC}}_W = {\mathsf{Park}}^{\text{NC}}_W(1)$ and 
${\mathsf{Park}}^{\text{alg}}_W = {\mathsf{Park}}^{\text{alg}}_W(1)$. 
There are simple formulas for the characters of these generalized parking spaces.
These character formulas lead to the following definition.

\begin{definition}
\label{def:phi}
Let $W$ be a finite complex reflection group acting on a complex vector space $V$.
For a positive integer $k$, 
let $\varphi^W_k$ and $\tilde{\varphi}^W_k$ be the class functions on $W$ defined by
\begin{gather}
\varphi^W_k(w) = k^{\dim V^w}
\quad(w \in W),
\label{eq:def_phi}
\\
\tilde{\varphi}^W_k(w) = \frac{ \det_V (1 - q^k w) }{ \det_V (1 - q w) }
\quad(w \in W),
\label{eq:def_phiq}
\end{gather}
where $V^w$ is the fixed-point subspace of $w$ and $q$ is an indeterminate.
We omit the superscript $W$ from $\varphi^W_k$ and $\tilde{\varphi}^W_k$ 
if there is no confusion.
\end{definition}

Note that $\tilde{\varphi}_k$ is a $q$-analogue of $\varphi_k$ 
in the sense that $\lim_{q \to 1} \tilde{\varphi}_k = \varphi_k$.

For a Weyl group $W$, Sommers \cite[Proposition~3.9]{Sommers} proved that, 
if $k$ is ``very good'', in particular if $k = h+1$, 
then $\varphi_k$ is the permutation character of $Q/kQ$.
For a real reflection group, it is shown (or conjectured) that
$\Comp[{\mathsf{Park}}^{\text{NC}}_W](p)$ and ${\mathsf{Park}}^{\text{alg}}_W(p)$ 
have the same character $\varphi_{ph+1}$ and 
the graded representation ${\mathsf{Park}}^{\text{alg}}_W(p)$ has 
the graded character $\tilde{\varphi}_{ph+1}$.
Here a graded representation of $W$ is a representation $U$ of $W$ equipped with 
a direct sum decomposition $U = \bigoplus_{i \in \Int} U_i$ into 
$W$-stable subspaces $U_i$, and the graded character of $U$ is defined as 
$\sum_{i \in \Int} \chi_i q^i$, where $\chi_i$ is the character of $U_i$.

However, $\varphi_k$ (or $\tilde{\varphi}_k$) does not necessarily come 
from genuine (graded) representations of $W$.
For example, if $W = \Sym_2$, then $\varphi_k$ is the character of 
a representation if and only if $k$ is odd.

\subsection{%
Main Theorem
}

The aim of this paper is to answer the following questions:

\begin{question}
\begin{enumerate}
\item[(1)]
When is $\varphi^W_k$ the character of a representation of $W$?
\item[(2)]
When is $\tilde{\varphi}^W_k$ is the graded character of 
a graded representation of $W$?
\end{enumerate}
\end{question}

We call a representation affording the character $\varphi_k$ or $\tilde{\varphi}_k$ 
a \emph{generalized parking space}.

To state our main theorem, we introduce generalizations of $W$-Catalan numbers.

\begin{definition}
\label{def:Catalan}
Let $W$ be a finite complex reflection group.
The \emph{generalized $q$-Catalan number} $\Cat_k(W,q)$ and 
the \emph{generalized dual $q$-Catalan number} $\Cat^*_k(W,q)$ 
are defined by
\begin{equation}
\label{eq:def_Catalan}
\Cat_k(W,q) = \prod_{i=1}^r \frac{ [k+d_i-1]_q }{ [d_i]_q },
\quad
\Cat^*_k(W,q) = q^N \prod_{i=1}^r \frac{ [k-d^*_i-1]_q }{ [d_i]_q },
\end{equation}
where $[m]_q = (1-q^m)/(1-q)$ is the $q$-integer, 
$d_1, \dots, d_r$ (resp. $d^*_1, \dots, d^*_r$) are the degrees (resp. the codegrees) 
of $W$ and $N$ is the number of reflecting hyperplanes.
\end{definition}

For example, $\Cat_{n+1}(\Sym_n,q)$ is MacMahon's $q$-analogue of the Catalan number, 
and $\Cat_{h+1}(W,q)$ is a $q$-analouge of $W$-Catalan number $\Cat_{h+1}(W,1)$.

The following is the main result of this paper.

\begin{theorem}
\label{thm:main}
Let $W$ be an irreducible finite complex reflection group and $k$ a positive integer.
Then we have
\begin{enumerate}
\item[(A)]
For any irreducible complex reflection group $W$, the following are equivalent:
\begin{enumerate}
\item[(i)]
$\tilde{\varphi}^W_k$ is the graded character of a graded representation of $W$.
\item[(ii)]
Both $\Cat_k(W,q)$ and $\Cat^*_k(W,q)$ are polynomials in $q$.
\item[(iii)]
$k$ satisfies the congruence condition given in Table~\ref{tab:condition}.
\begin{table}[htbp]
\caption{Condition in Theorem~\ref{thm:main} (iii)}
\label{tab:condition}
$$
\begin{array}{c|c}
\text{group} & \text{condition on $k$} \\
\hline
\Sym_n & \gcd(k,n) = 1 \\
\hline
\begin{array}{c} G(m,p,n) \\ (n \ge 3 \text{ or } p<m) \end{array} & k \equiv 1 \bmod m \\
\hline
G(m,m,2) & k \equiv \pm 1 \bmod m \\
\hline
C_m & k \equiv 1 \bmod m \\
\hline
G_4 & k \equiv 1, 3 \bmod 6 \\
G_5 & k \equiv 1 \bmod 6 \\
G_6 & k \equiv 1, 9 \bmod 12 \\
G_7 & k \equiv 1 \bmod 12 \\
G_8 & k \equiv 1, 5 \bmod 12 \\
G_9 & k \equiv 1, 17 \bmod 24 \\
G_{10} & k \equiv 1 \bmod 12 \\
G_{11} & k \equiv 1 \bmod 24 \\
G_{12} & k \equiv 1, 11, 17, 19 \bmod 24 \\
G_{13} & k \equiv 1, 17 \bmod 24 \\
G_{14} & k \equiv 1, 19 \bmod 24 \\
G_{15} & k \equiv 1 \bmod 24 \\
G_{16} & k \equiv 1, 11 \bmod 30 \\
G_{17} & k \equiv 1, 41 \bmod 60 \\
G_{18} & k \equiv 1 \bmod 30 \\
G_{19} & k \equiv 1 \bmod 60 \\
G_{20} & k \equiv 1, 19 \bmod 30 \\
G_{21} & k \equiv 1, 49 \bmod 60 \\
G_{22} & k \equiv 1, 29, 41, 49 \bmod 60 \\
\hline
G_{23} = W(H_3) & k \equiv 1, 5, 9 \bmod 10 \\
G_{24} & k \equiv 1, 9, 11 \bmod 14 \\
G_{25} & k \equiv 1 \bmod 6 \\
G_{26} & k \equiv 1 \bmod 6 \\
G_{27} & k \equiv 1, 19, 25 \bmod 30 \\
\hline
G_{28} = W(F_4) & k \equiv 1, 5 \bmod 6 \\
G_{29} & k \equiv 1, 9, 13, 17 \bmod 20 \\
G_{30} = W(H_4) & k \equiv 1, 11, 19, 29 \bmod 30 \\
G_{31} &  k \equiv 1, 13 ,17, 29, 37, 41, 49, 53 \bmod 60 \\
G_{32} & k \equiv 1, 7, 13, 19 \bmod 30 \\
\hline
G_{33} & k \equiv 1 \bmod 6 \\
\hline
G_{34} & k \equiv 1, 13, 19, 25, 31, 37 \bmod 42 \\
G_{35} = W(E_6) & k \equiv 1, 5 \bmod 6 \\
\hline
G_{36} = W(E_7) & k \equiv 1, 5 \bmod 6 \\
\hline
G_{37} = W(E_8) & k \equiv 1, 7, 11, 13, 17, 19, 23, 29 \bmod 30
\end{array}
$$
\end{table}
\end{enumerate}
\item[(B)]
Except for dihedral groups $W = G(m,m,2) = D_{2m}$, 
the conditions in (A) are equivalent to the following:
\begin{enumerate}
\item[(iv)]
$\varphi^W_k$ is the character of a representation of $W$.
\item[(v)]
$\varphi^W_k$ is the character of a permutation representation of $W$.
\end{enumerate}
\item[(C)]
If $W = G(m,m,2) = D_{2m}$, then the following conditions are equivalent:
\begin{enumerate}
\item[(iii')]
$k = 1$, or $k \ge m-1$ and $k$ satisfies the congruence
$$
k^2 \equiv 1 \ \begin{cases}
 \bmod \ 2m &\text{if $m$ is even,} \\
 \bmod \ m &\text{if $m$ is odd.}
\end{cases}
$$
\item[(iv)]
$\varphi^W_k$ is the character of a representation of $W$.
\item[(v)]
$\varphi^W_k$ is the character of a permutation representation of $W$.
\end{enumerate}
\end{enumerate}
\end{theorem}

We should remark that $\varphi_k$ does not come from a genuine representation of $W$ 
even if $k$ is relatively prime to the Coxeter number $h$.
For example, if $W$ is the Coxeter group of type $H_3$, 
then $k = 3$ is relatively prime to the Coxeter number $h = 10$, 
but $\varphi_3$ is not the character of a representation of $W$.

Some partial results on the existence of generalized parking spaces can be found in the literatures.
Haiman \cite[Propositio~2.4.1]{Haiman} proved the equivalence of (i), (ii), and (iii) 
for the symmetric groups.
Sommers \cite[Proposition~3.9]{Sommers} proved the implication (iii) $\implies$ (v) for Weyl groups.
Note that $k$ is ``very good'' in the sense of \cite{Sommers} 
if and only if $k$ satisfies the condition in (iii).
Also it is known that some of $\tilde{\varphi}_k$, up to a power of $q$, 
appear as the characters of finite dimensional representations of rational Cherednik algebras.
See \cite{BEG} for example.

\subsection{%
Strategy of proof
}

Our proof of Theorem~\ref{thm:main} proceeds as follows.
We prove the following six implications
\begin{enumerate}
\item[(a)] (i) $\implies$ (ii),
\item[(b)] (ii) $\implies$ (iii),
\item[(c)] (iii) $\implies$ (i),
\item[(d)] (iii) (or (iii') for dihedral groups) $\implies$ (v) ,
\item[(e)] (v) $\implies$ (iv),
\item[(f)] (iv) $\implies$ (iii) (or (iii') for dihedral groups).
\end{enumerate}
Among these implications, (e) is obvious, 
and (a) follows from the fact that the generalized $q$-Catalan numbers 
$\Cat_k(W,q)$ and $\Cat^*_k(W,q)$ appear 
as multiplicities of some irreducible characters in $\tilde{\varphi}_k$
(See Proposition~\ref{prop:Catalan=mult}).
The other implications will be proved in a case-by-case manner.

The implication (b) is proved by decomposing $\Cat_k(W,q)$ and $\Cat^*_k(W,q)$ 
into the cyclotomic polynomials.

For the implications (c) and (f), we write
$$
\varphi_k = \sum_{\chi \in \Irr(W)} m_k^\chi \chi,
\quad
\tilde{\varphi}_k = \sum_{\chi \in \Irr(W)} \tilde{m}_k^\chi \chi,
$$
where $\Irr(W)$ is the set of all irreducible characters of $W$.
Then the condition (iv) (resp. (i)) is equivalent to saying that 
$m_k^\chi$ is a nonnegative integers 
(resp. $\tilde{m}_k^\chi$ is a polynomials in $q$ with nonnegative integer coefficients) 
for all $\chi \in \Irr(W)$.
And the multiplicities $m_k^\chi$ and $\tilde{m}_k^\chi$ can be explicitly computed.
If $W = \Sym_n$ or $G(m,p,n)$, then these multiplicities are given in terms of 
specializations of Schur functions (see Proposition~\ref{prop:decomp_S} 
and Corollary~\ref{cor:decomp_G2}).
For the exceptional groups, we can use the program GAP \cite{GAP3} together with CHEVIE \cite{CHEVIE} 
to compute these multiplicities.
By using these explicit formulas for the multiplicities, we can prove (c) and (f).

In order to prove (d), we appeal to a result of Orlik--Solomon \cite{OS1, OS2}, 
which provides an expression of $\varphi_k$ as a linear combination of 
certain permutation characters $\eta_j$:
$$
\varphi_k
 = 
\sum_{j=1}^s n_j(k) \eta_j,
$$
where $n_j(k)$ is a polynomial in $k$.
(See Proposition~\ref{prop:OS}.)
Thus we can prove (d) by showing that all $n_j(k)$'s are nonnegative integers 
if $k$ satisfies the condition in (iii).

In the proof for the symmetric groups, a key role is played by the following results 
on the greatest common divisors of specialized Schur functions, 
which is interesting in itself.

\begin{theorem}
(Theorem~\ref{thm:gcd} below)
Let $k$ and $n$ be positive integers.
Then we have
\begin{gather*}
\gcd{}_\Int \{ s_\lambda (\underbrace{1, \dots, 1}_k) : \lambda \vdash n \}
 =
\frac{k}{\gcd(n,k)},
\\
\gcd{}_{\Rat[q]} \{ s_\lambda (1, q, \dots, q^{k-1}) : \lambda \vdash n \}
 =
\frac{[k]_q}{[\gcd(n,k)]_q}.
\end{gather*}
\end{theorem}

\subsection{%
Organization
}

The remaining of this paper is organized as follows.
In Section~2, we prepare several general results which are useful in the proof 
of our main result (Theorem~\ref{thm:main}), especially for exceptional groups.
Sections~3, 4,  5, and 6 are devoted to the proof of the main result 
in the cases where 
$W = \Sym_n$, $G(m,p,n)$, $D_{2m}$ and exceptional groups, respectively.
In Subsection~3.2, we prove the formulae for the greatest common divisors of 
specialized Schur functions.
In Section~7, we conclude with some problems.

\subsection{%
Notations
}

In this paper, $\Nat$ and $\Nat[q]$ denote the set of nonnegative integers 
and the set of polynomials in $q$ with nonnegative integer coefficients respectively.
\section{%
Preliminaries
}

In this section, we review two general results for complex reflection groups 
and present several criteria for nonnegativity and polynomiality.

\subsection{%
$q$-Catalan numbers and multiplicities
}

We refer the readers to \cite{LT} for the theory of complex reflection groups.

Let $V$ be a complex vector space of dimension $r$.
For an linear map $w : V \to V$, we put
$$
V^w = \{ v \in V : w(v) = v \}.
$$
A complex reflection on $V$ is an invertible linear map $w : V \to V$ 
satisfying $\dim V^w = \dim V - 1$.
Such a hyperplane $V^w$ is called a reflecting hyperplane.
A \emph{complex reflection group} is a subgroup $W$ of $\GL(V)$ which is 
generated by complex reflections.
Then $V$ is called the reflection representation of $W$.
We say that $W$ is \emph{irreducible} if $V$ has no $W$-stable subspace 
other than $\{ 0 \}$ and $V$.
Shephard--Todd \cite{ST} classified the irreducible finite complex reflection groups 
into three infinite families
\begin{enumerate}
\item[1.]
the symmetric groups $\Sym_n$,
\item[2.]
the groups $G(m,p,n)$,
\item[3.]
the cyclic groups $C_m$,
\end{enumerate}
and 34 exceptional groups $G_4, \dots, G_{37}$.

Let $W$ be a complex reflection group with $r$-dimensional reflection representation $V$.
Given an $n$-dimensional $W$-module $M$, 
the $M$-exponents $m_1(M) \le m_2(M) \le \dots \le m_n(M)$ are defined by
$$
\sum_{i=1}^n q^{m_i(M)}
 =
\sum_{d \ge 0} \dim \Hom_W(\mathcal{H}_d, M) q^d,
$$
where $\mathcal{H} = \bigoplus_{d \ge 0} \mathcal{H}_d$ is the space of $W$-harmonic polynomials.
Then the \emph{degrees} $(d_1, \dots, d_r)$ and the \emph{codegrees} $(d^*_1, \dots, d^*_r)$ 
of $W$ are given in terms of the $M$-exponents for $M = V$ and $V^*$:
$$
d_i = m_i(V) + 1,
\quad
d^*_i = m_i(V^*) - 1
\quad
(1 \le i \le r).
$$
See \cite[pp.~274--275]{LT} for the data of degrees and codegrees.
It is known that the number of reflecting hyperplanes are given by
\begin{equation}
\label{eq:no_refl_hyperplane}
N = \sum_{i=1}^r (d^*_i+1).
\end{equation}

It is convenient to work with the following generalization of $\tilde{\varphi}_k$.

\begin{definition}
\label{def:phiqu}
For indeterminates $q$ and $u$, let $\hat{\varphi}$ be the class function on $W$ defined by
\begin{equation}
\label{eq:def_phiqu}
\hat{\varphi}(w)
 =
\frac{ \det_V (1 - u w) }
     { \det_V (1 - q w) }
\quad(w \in W).
\end{equation}
\end{definition}

Note that $\tilde{\varphi}_k$ is obtained from $\hat{\varphi}$ 
by substituting $u = q^k$.
If we denote by $\Irr(W)$ the set of irreducible characters of $W$, then we can write
\begin{equation}
\label{eq:irr_decomp}
\hat{\varphi}
 =
\sum_{\chi \in \Irr(W)} \hat{m}^\chi \chi,
\quad
\tilde{\varphi}_k
 =
\sum_{\chi \in \Irr(W)} \tilde{m}_k^\chi \chi,
\quad
\varphi_k
 =
\sum_{\chi \in \Irr(W)} m_k^\chi \chi,
\end{equation}
where $\hat{m}^\chi \in \Comp(q)[u]$, $\tilde{m}_k^\chi \in \Comp(q)$, 
and $m_k^\chi \in \Comp$.
It follows from 
$$
\hat{m}^\chi
 =
\frac{  1}{ \# W }
\sum_{w \in W}
 \hat{\varphi}(w) \overline{\chi(w)},
$$
that $\hat{m}^\chi$ is a polynomial in $u$ of degree $\le r$.
Similarly we see that $m^\chi_k$ is a polynomial in $k$ of degree $\le k$.

\begin{prop}
\label{prop:Catalan=mult}
The multiplicities of the trivial character $\triv$ 
and the determinant character $\det = \det_V$ in 
$\hat{\varphi}$ are given by
\begin{equation}
\label{eq:mult_phiqu}
\hat{m}^{\triv}
 =
\prod_{i=1}^r
 \frac{ 1 - u q^{d_i - 1} }
      { 1 - q^{d_i} },
\quad
\hat{m}^{\det}
 =
\prod_{i=1}^r
 \frac{ 1 - u q^{- d_i^* - 1} }
      { 1 - q^{d_i} }.
\end{equation}
In particular, we have
\begin{equation}
\label{eq:Catalan=mult}
\tilde{m}_k^{\triv}
 =
\Cat_k(W,q),
\quad
\tilde{m}_k^{\det}
 =
\Cat^*_k(W,q).
\end{equation}
\end{prop}

\begin{demo}{Proof}
For a linear character $\chi$ of $W$, we have
(\cite[(2.3)]{Lehrer})
$$
\frac{1}{\# W}
\sum_{w \in W}
 \hat{\varphi}(w) \chi(w)
 =
q^{m(\Comp_\chi)}
\prod_{i=1}^r
 \frac{ 1 - u q^{m_i(V \otimes \Comp_\chi) - m(\Comp_\chi)} }
      { 1 - q^{d_i} },
$$
where $\Comp_\chi$ is the one-dimensional $W$-module affording $\chi$, 
and $m(\Comp_\chi)$, $\{ m_i (V \otimes \Comp_\chi) \}_{i=1}^r$ denote the exponents 
of $\Comp_\chi$ and $V \otimes \Comp_\chi$ respectively.
We apply this formula to $\chi = \triv$ and $\det^{-1}$.
By using the relations (see \cite{Lehrer} for example)
\begin{gather*}
m(\Comp_{\triv}) = 0,
\quad
m(\Comp_{\det^{-1}}) = N,
\\
m_i(V \otimes \Comp_{\triv}) = m_i(V),
\quad
m_i(V \otimes \Comp_{\det^{-1}}) = N - m_{r+1-i}(V^*),
\end{gather*}
we obtain the formulas (\ref{eq:mult_phiqu}).
\end{demo}

\subsection{%
Orlik--Solomon formula
}

Let $W$ be a complex reflection group acting on $V$, 
and $\mathcal{A}$ the set of reflecting hyperplanes.
We denote by $L(\mathcal{A})$ the intersection lattice of $\mathcal{A}$.
That is, $L(\mathcal{A})$ is the set of all intersections of hyperplanes 
in $\mathcal{A}$ ordered by reverse inclusion.
Then $W$ acts on $L(\mathcal{A})$.

\begin{prop}
\label{prop:OS}
(Orlik--Solomon \cite{OS1, OS2})
Let $X_1, \dots, X_s$ be a complete set of representatives of $W$-orbits on $L(\mathcal{A})$.
For $1 \le j \le s$, we put
\begin{gather*}
W_j = \{ w \in W : X_j \subset V^w \},
\\
\mathcal{A}_j = \{ H \cap X_j : H \in \mathcal{A}, \, H \not\supset X_j \},
\end{gather*}
and denote by $\eta_j$ the permutation character of $W$ on $W/W_j$.
Then we have
\begin{equation}
\label{eq:OS}
\varphi_k
 = 
\sum_{j=1}^s
 \frac{ \chi(\mathcal{A}_j, k) }{ [N_W(W_j) : W_j] } \eta_j,
\end{equation}
where $\chi(\mathcal{A}_j, t)$ is the characteristic polynomial of 
the hyperplane arrangement $\mathcal{A}_j$, and 
$N_W(W_j)$ is the normalizer of $W_j$ in $W$.
\end{prop}

The polynomials $\chi(\mathcal{A}_j,t)$ are monic polynomials and 
the roots are known to be positive integers.
For exceptional groups $W = G_{23}, \dots, G_{37}$ of rank $\ge 3$, 
the $W$-orbits on $L(\mathcal{A})$ and the roots of the corresponding 
characteristic polynomials $\chi(\mathcal{A}_j,t)$ are listed in 
\cite[Tables~III--VIII]{OS1} (Coxeter groups of type $H_3, H_4, F_4, E_6, E_7, E_8$) 
and \cite[Tables~3--11]{OS2} (other groups).
For the remaining rank $2$ groups $W = G_4, \dots, G_{22}$, 
one can find the lists of the orbits $\mathcal{O}_j$ with $W_j$ and their sizes $\# \mathcal{O}_j$ 
in \cite[Table~2]{OS2}, and the index is computed by
$$
[N_W(W_j) : W_j]
 = 
\frac{ \# W }
     { \# W_j \cdot \# \mathcal{O}_j }.
$$
For the groups $\Sym_n$ and $G(m,1,n)$, we can use the theory of 
symmetric functions to prove the corresponding results.
See Propositions~\ref{prop:decomp_S} (2) and \ref{prop:decomp_G3}.

\subsection{%
Nonnegativity
}

In this subsection, we present three general lemmas.
The proofs are easy, so we omit them.

In order to prove the implication (iv) $\implies$ (iii) of Theorem~\ref{thm:main} 
for exceptional groups, we need to find all integers $k$ satisfying $m^\chi_k \in \Nat$ 
for a fixed irreducible character $\chi$ of $W$.
Since $m^\chi_k$ is a polynomial in $k$, we can appeal to the following lemma to 
reduce the proof to a finite amount of computation.

\begin{lemma}
\label{lem:period}
Let $f \in \Comp[t]$ be a polynomial in $t$.
Assume that there is a positive integer $L$ such that $f(t+L) - f(t)$ maps 
nonnegative integers to nonnegative integers.
Then, for an integer $i$ with $1 \le i \le L$, the following are equivalent:
\begin{enumerate}
\item[(i)]
$f(i) \in \Nat$.
\item[(ii)]
$f(k) \in \Nat$ for any integer $k$ with $k \equiv i \bmod L$.
\end{enumerate}
\end{lemma}

The assumption in the lemma above can be checked by using the following lemma.

\begin{lemma}
\label{lem:NtoN}
Let $g \in \Comp[t]$ be a polynomial in $t$ of degree $\le r$.
We define $b_0, b_1, \dots, b_r$ recursively by the relations
$$
b_0 = g(0),
\quad
b_i = g(i) - \sum_{j=0}^{i-1} b_j \binom{i}{j}.
$$
Then we have
$$
g(t) = \sum_{i=0}^r b_i \binom{t}{i},
$$
and, if $b_0, \dots, b_r \in \Nat$, then 
$g(t)$ maps nonnegative integers to nonnegative integers.
\end{lemma}

This lemma can also be used in the proof of (iii) $\implies$ (v), 
where we prove that $\chi(\mathcal{A}_j,pH+i)/[N_W(W_i):W_i] \in \Nat$ 
for all $p \in \Nat$.

The following lemma is a $q$-analogue of Lemma~\ref{lem:NtoN}.
This lemma will be used in the proof of (iii) $\implies$ (i), 
where we show that $\tilde{m}^\chi_{pH+i} = \hat{m}^\chi(q^{pH+i}) \in \Nat[q]$ 
for all $p \in \Nat$.

\begin{lemma}
\label{lem:NtoNq}
Let $h(u)$ be a polynomial in $u$ of degree $\le r$ with coefficients in $\Comp(q)$, 
and $M$ a positive integer.
Define $c_0, c_1, \dots, c_r$ recursively by the relations
$$
c_0 = h(1),
\quad
c_i = h(q^{iM}) - \sum_{j=0}^{i-1} \gauss{i}{j}{q^M} c_j,
$$
where $\gauss{i}{j}{q^M}$ is the $q$-binomial coefficient in base $q^M$, i.e.,
$$
\gauss{i}{j}{q^M}
 =
\prod_{l=1}^j
 \frac{ 1 - q^{(i-l+1)M} }
      { 1 - q^{lM} }.
$$
Then we have
$$
h(u)
=
\sum_{i=0}^r c_i 
 \prod_{l=1}^{i} 
  \frac{ 1 - u q^{(-l+1)M} }
       { 1 - q^{lM} },
$$
and, if $c_0, c_1, \cdots, c_r \in \Nat[q]$, 
then $h(q^{pM}) \in \Nat[q]$ for all nonnegative integers $p$.
\end{lemma}

\subsection{%
Polynomiality of $q$-Catalan numbers
}

In the proof of (ii) $\implies$ (iii) of Theorem~\ref{thm:main}, 
we need to find all integers $k$ such that $\Cat_k(W,q)$ and $\Cat^*_k(W,q)$ 
are both polynomials in $q$.
By using (\ref{eq:no_refl_hyperplane}), we have
$$
\Cat^*_k(W,q)
 = 
\prod_{i=1}^r 
 \frac{ q^{d^*_i+1} - q^k }
      { 1 - q^{d_i} }.
$$
Thus $\Cat^*_k(W,q)$ is a polynomial in $q$ if and only if 
$\Cat^*_k(W,q)$ is a Laurent polynomial in $q$.

Let $\Phi_d(q)$ denotes the $d$-th cyclotomic polynomial in $q$. 
Then $\Phi_d(q) \in \Int[q]$ is an irreducible polynomial and we have
\begin{equation}
1 - q^a = \prod_{d \mid a} \Phi_d(q),
\label{eq:cyclotomic}
\end{equation}
where $d$ runs over all divisors of a positive integer $a$.
The Laurent polynomiality of $\Cat_k(W,q)$ and $\Cat^*_k(W,q)$ can be checked 
by using the following lemma:

\begin{lemma}
\label{lem:pol_Catalan}
Given integers $a_1, \cdots, a_r$, and positive integers $b_1, \cdots, b_r$,
we consider the rational function given by
$$
c(q)
 =
\prod_{i=1}^r \frac{ [a_i]_q }{ [b_i]_q }.
$$
For a positive integer $d$, we put
$$
N(d) = \# \{ i : d \mid a_i \},
\quad
D(d) = \# \{ i : d \mid b_i \}.
$$
Then $c(q)$ is a Laurent polynomial in $q$ if and only if 
some of $a_i$'s are 0 or $N(d) \ge D(d)$ for all divisors $d \in T$, where
$$
T = \bigcup_{i=1}^s \{ d : d \mid b_i \}.
$$
\end{lemma}

\begin{demo}{Proof}
Note that $c(q) = 0$ if and only if some of $a_i$'s are $0$.
So we may assume that all $a_i$'s are nonzero.
By using the relation $[-a]_q = -q^{-a} [a]_q$ and (\ref{eq:cyclotomic}), we have
$$
c(q)
=
\ep q^{-S} \prod_d \Phi_d(q)^{N(d) - D(d)},
$$
where $d$ runs over all positive integers, and 
$\ep \in \{ 1, -1 \}$, $S = \sum_{i : a_i < 0} a_i$.
Since $D(d) = 0$ for $d \not\in T$, we see that 
$c(q)$ is a Laurent polynomial if and only if $N(d)-D(d) \ge 0$ for all $d \in T$.
\end{demo}

\section{%
Symmetric groups
}

In this section, we prove the main result (Theorem~\ref{thm:main}) 
in the case where $W = \Sym_n$ is the symmetric group. 
As a key to our proof, we give formulae for the greatest common divisors of 
specialized Schur functions.

\subsection{%
Representations of $\Sym_n$ and symmetric functions
}

In this subsection, we review several facts concerning representations 
of $\Sym_n$ and symmetric functions.
See \cite[Chapter~1]{Macdonald} or \cite[Chapter~7]{EC2} for details.

The symmetric group $\Sym_n$ of $n$ letters is realized as an irreducible 
complex reflection group acting on $V = \{ x \in \Comp^n : x_1 + \dots + x_n = 0 \}$.

A partition of a nonnegative integer $n$ is a weakly decreasing sequence 
$\lambda = (\lambda_1, \lambda_2, \dots)$ of nonnegative integers satisfying 
$\sum_{i \ge 1} \lambda_i = n$.
Then we write $|\lambda| = n$ and $\lambda \vdash n$.
The length of a partition $\lambda$, denoted by $l(\lambda)$, is the number 
of nonzero entries of $\lambda$.
Let $\Par_n$ be the set of partitions of $n$.
We often identify a partition $\lambda$ with the finite sequence 
$(\lambda_1, \dots, \lambda_k)$ with $k \ge l(\lambda)$.

The conjugacy classes of $\Sym_n$ are parametrized by partitions of $n$, 
called the \emph{cycle type}.
A permutation $w \in \Sym_n$ has the cycle type $\mu \vdash n$, denoted by $\type(w)$, 
if $w$ is decomposed into the product of disjoint cycles of lengths $\mu_1, \mu_2, \dots$. 
Then we have
$$
\det{}_V (1 - tw )
 = 
\frac{ \prod_{i \ge 1} (1 - t^{\mu_i}) }
     { 1-t },
$$
and
\begin{equation}
\label{eq:value_S}
\tilde{\varphi}_k(w)
 =
\frac{ \prod_{i \ge 1} [k]_{q^{\mu_i}} }
     { [k]_q },
\quad
\varphi_k(w)
 =
k^{l(\mu) - 1}.
\end{equation}

The irreducible characters of $\Sym_n$ are also indexed by partitions of $n$.
Let $\chi^\lambda$ be the irreducible character of $\Sym_n$ corresponding to 
a partition $\lambda$, 
and denote by $\chi^\lambda_\mu$ the character value at an element of cycle type $\mu$.

Let $R(\Sym_n)$ be the vector space of complex-valued class functions on $\Sym_n$.
Then the direct sum $R(\Sym_\bullet) = \bigoplus_{n \ge 0} R(\Sym_n)$ 
has a commutative associative graded algebra structure 
with respect to the product defined by
$$
f \cdot g = \Ind_{\Sym_n \times \Sym_l}^{\Sym_{n+l}} (f \times g),
$$
where $f \in R(\Sym_n)$ and $g \in R(\Sym_l)$.

Let $\Lambda$ be the ring of symmetric functions in infinitely many variables 
$\vectx = (x_1, x_2, \dots)$ with complex coefficients. 
We follow \cite{Macdonald} for the notations of symmetric functions.
For example, we denote by $s_\lambda$ the Schur function associated to a partition $\lambda$.
The following identity is known as the Cauchy identity:

\begin{prop}
\label{prop:Cauchy}
For two sets of variables $\vectx$ and $\vecty$, we have
\begin{equation}
\label{eq:Cauchy}
\prod_{i,j} ( 1 - x_i y_j )^{-1}
 =
\sum_\lambda s_\lambda(\vectx) s_\lambda(\vecty)
 =
\sum_\lambda h_\lambda(\vectx) m_\lambda(\vecty),
\end{equation}
where $h_\lambda$ and $m_\lambda$ are the complete and monomial symmetric functions respectively.
\end{prop}

Define a linear map $\ch : R(\Sym_\bullet) \to \Lambda$, called the \emph{Frobenius characteristic}, by
$$
\ch(f) = \frac{1}{n!} \sum_{w \in \Sym_n} f(w) p_{\type(w)}
$$
for $f \in R(\Sym_n)$.
Then we have

\begin{prop}
\label{prop:Frobenius_S}
\begin{enumerate}
\item[(1)]
The Frobenius characteristic $\ch$ is an algebra isomorphism.
\item[(2)]
For a partition $\lambda$ of $n$, 
the images of the irreducible character $\chi^\lambda$ and 
the permutation character $\eta^\lambda$ on $\Sym_n/\Sym_\lambda$ are given by
\begin{equation}
\label{eq:Frobenius_S1}
\ch(\chi^\lambda) = s_\lambda,
\quad
\ch(\eta^\lambda) = h_\lambda.
\end{equation}
\item[(3)]
For a partition $\mu$ of $n$, we have
\begin{equation}
\label{eq:Frobenius_S2}
p_\mu (\vectx)
 = 
\sum_{\lambda \vdash n} \chi^\lambda_\mu s_\lambda(\vectx),
\end{equation}
where $\lambda$ runs over all partitions of $n$.
\end{enumerate}
\end{prop}

\begin{definition}
\label{def:specialization}
For a symmetric function $f \in \Lambda$ and a positive integer $k$, 
we denote by $f(1^k)$ the specialization of $f$ with $x_1 = \dots = x_k = 1$ 
and $x_{k+1} = x_{k+2} = \dots = 0$.
And we denote by $f(1, q, \dots, q^{k-1})$ the principal specialization of $f$ 
with $x_i = q^{i-1}$ ($1 \le i \le k$ and $x_i = 0$ ($i > k$).
\end{definition}

By using the Frobenius formula (\ref{eq:Frobenius_S2}) and the Cauchy identity (\ref{eq:Cauchy}), 
we obtain the following expressions of $\varphi_k$ and $\tilde{\varphi}_k$.

\begin{prop}
(See \cite{Stanley1997})
\label{prop:decomp_S}
\begin{enumerate}
\item[(1)]
The class function $\tilde{\varphi}_k$ is expressed in terms of irreducible characters as
\begin{equation}
\label{eq:decomp_S1}
\tilde{\varphi}_k
 =
\sum_{\lambda \vdash n}
 \frac{ s_\lambda(1, q, \dots, q^{k-1}) }{ [k]_q }
 \chi^\lambda.
\end{equation}
In particular, we have
\begin{equation}
\label{eq:decomp_S2}
\varphi_k
 = 
\sum_{\lambda \vdash n}
 \frac{ s_\lambda(1^k) }{ k }
 \chi^\lambda.
\end{equation}
\item[(2)]
The class function $\varphi_k$ is expressed in terms of permutation characters as
\begin{equation}
\label{eq:decomp_S3}
\varphi_k
 =
\sum_{\lambda \vdash n}
 \frac{ m_\lambda(1^k) }{ k }
 \eta_\lambda.
\end{equation}
\end{enumerate}
\end{prop}

\begin{demo}{Proof}
(1)
By substituting $x_i = q^{i-1}$ ($1 \le i \le k$) and $x_i = 0$ ($i \ge k+1$) 
in the Frobenius formula (\ref{eq:Frobenius_S2}) and comparing with (\ref{eq:value_S}), 
we obtain (\ref{eq:decomp_S1}).

(2)
By substituting $y_1 = \dots = y_k = 1$ and $y_{k+1} = \dots = 0$ in the Cauchy identity 
(\ref{eq:Cauchy}), and by using (\ref{eq:Frobenius_S1}), we have
$$
\ch (\varphi_k)
 =
\sum_{\lambda \vdash n} \frac{ s_\lambda(1^k) }{k} s_\lambda
 =
\sum_{\lambda \vdash n} \frac{ m_\lambda(1^k) }{k} h_\lambda
 =
\sum_{\lambda \vdash n} \frac{ m_\lambda(1^k) }{k} \ch (\eta_\lambda).
$$
Since $\ch$ is an isomorphism, we obtain (\ref{eq:decomp_S3}).
\end{demo}

\subsection{%
Greatest common divisor of specialized Schur functions
}

It follows from Proposition~\ref{prop:decomp_S} (1) that 
$\varphi_k$ is the character of a representation of $\Sym_n$ 
if and only if 
$s_\lambda(1^k)$ is divisible by $k$ for any partitions $\lambda \vdash n$.
Similarly, $\tilde{\varphi}_k$ is the graded character of a graded representation 
of $\Sym_n$ if and only if $s_\lambda(1, q, \dots, q^{k-1})$ is divisible by $[k]_q$ 
and the quotient is a polynomial with nonnegative integer coefficients 
for any $\lambda \vdash n$.
So the following theorem plays an indispensable role in the proof of Theorem~\ref{thm:main}.

\begin{theorem}
\label{thm:gcd}
Let $k$ and $n$ be positive integers.
\begin{enumerate}
\item[(1)]
In the integer ring $\Int$, we have
\begin{equation}
\label{eq:gcd1}
\gcd{}_\Int \{ s_\lambda (1^k) : \lambda \vdash n \}
 =
\frac{k}{\gcd(n,k)}.
\end{equation}
\item[(2)]
In the polynomial ring $\Rat[q]$, we have
\begin{equation}
\label{eq:gcd2}
\gcd{}_{\Rat[q]} \{ s_\lambda (1, q, \dots, q^{k-1}) : \lambda \vdash n \}
 =
\frac{[k]_q}{[\gcd(n,k)]_q}.
\end{equation}
Here the greatest common divisor is taken to be a monic polynomial.
\end{enumerate}
\end{theorem}

\begin{remark}
It should be noted that (1) is not an immediate consequence of its ``q-analogue'' (2).
For example, if $f(q) = (q^2 + 1)(q + 1)^2$ and $g(q) = (q + 1)^3$, then 
we have $h(q) = \gcd_{\Rat[q]} \{ f(q), g(q) \} = (q+1)^2$ and $h(1) = 4$, 
while $\gcd_{\Int} \{ f(1), g(1) \} = 8$.
\end{remark}

In the proof of Theorem~\ref{thm:gcd}, we will use the following two lemmas.

\begin{lemma}
\label{lem:root}
(Haiman \cite[Proof of Proposition~2.5.1]{Haiman})
Let $d$, $k$, $r$ be positive integers such that $d$ divides $k$.
\begin{enumerate}
\item[(1)]
If $d \mid r$, then no primitive $d$-th root of $1$ is 
a root of $h_r (1, q, \dots, q^{k-1})$.
\item[(2)]
If $d \nmid r$, then any primitive $d$-th root of $1$ is 
a simple root of $h_r (1, q, \dots, q^{k-1})$.
\end{enumerate}
\end{lemma}

\begin{demo}{Proof}
If $r = sd+t$ with $s \in \Nat$ and $0 \le t < d$, then we have
$$
h_r(1, q, \dots, q^{k-1})
 =
\prod_{i=0}^{s-1} \prod_{j=1}^d 
 \frac{ 1 - q^{id+j+k-1} }
      { 1 - q^{id+j} }
\cdot
\prod_{j=1}^t
 \frac{ 1 - q^{sd+j+k-1} }
      { 1 - q^{sd+j} },
$$
and the proof follows from this expression.
\end{demo}

For a prime $p$ and an integer $x$, 
we denote by $\nu_p(x)$ the highest power of $p$ dividing $x$.

\begin{lemma}
\label{lem:Kummer}
(Kummer \cite[p.~116]{Kummer})
Let $p$ be a prime, and $m$, $r$ positive integers with $m \ge r$.
Then $\nu_p \left( \binom{m}{r} \right)$ is equal to the number of borrows 
required when subtracting $r$ from $m$ in the base $p$ representation.
\end{lemma}

\begin{demo}{Proof of Theorem~\ref{thm:gcd}}
In the proof, we put $d = \gcd (n,k)$.

First we prove (2).
Let $g(q)$ be the greatest common divisor of $s_\lambda(1, q, \dots, q^{k-1})$'s with $\lambda \vdash n$.
Since $\{ h_\lambda : \lambda \vdash n \}$ is another $\Int$-basis of 
$\Lambda_\Int = \sum_{\lambda \vdash n} \Int s_\lambda$, 
we have $g(q)$ is the greatest common divisor of $h_\lambda(1, q \dots, q^{k-1})$'s.
Also we have
$$
\frac{ [k]_q }{ [d]_q }
 =
\prod_{d \mid k, \, d \nmid n} \Phi_d(q)
 =
\prod_{d \mid k, \, d \nmid n} \prod_\zeta (q - \zeta),
$$
where $\Phi_d$ is the $d$-th cyclotomic polynomial 
and $\zeta$ runs over all primitive $d$-th roots of unity.
In order to prove $g(q) = [k]_q/[d]_q$, it is enough to show the following two claims:

\begin{description}
\item[Claim~1]
We have
\begin{multline*}
\{ z \in \Comp : \text{$z$ is a common root of $h_\lambda (1, q, \dots, q^{k-1})
 \; (\lambda \vdash n)$} \}
 \\
=
\bigsqcup_{d \mid k, \, d \nmid n}
\{ z \in \mathbb{C} : \text{$z$ is a primitive $d$-th root of $1$} \}.
\end{multline*}
\item[Claim~2]
If $z$ is a common root of $h_\lambda (1, q, \dots, q^{k-1})$ ($\lambda \vdash n$),
then $z$ is a simple root of $h_\mu (1, q, \dots, q^{k-1})$ for some $\mu \vdash n$.
\end{description}

\begin{demo}{Proof of Claim~1}
We put
\begin{align*}
C
 &=
\{ z \in \Comp :
 \text{$z$ is a common root of $h_\lambda (1, q, \dots, q^{k-1}) \; (\lambda \vdash n)$}
\},
\\
D
 &=
\bigsqcup_{d \mid k \, d \nmid n}
\{ z \in \mathbb{C} \mid \text{$z$ is a primitive $d$-th root of $1$} \}.
\end{align*}
Let $z \in C$.
Since $z$ is a root of $h_{(1^n)} (1, q, \dots, q^{k-1})$ and
$$
h_{(1^n)} (1, q, \dots, q^{k-1})
 =
\left( h_1 (1, q, \dots, q^{k-1}) \right)^n
 =
\left( \prod_{d \mid k, \, d \neq 1} \Phi_d (q) \right)^n,
$$
we see that $z$ is a primitive $d$-th root of $1$ for some $d$ dividing $k$.
Since $z$ is also a root of $h_n(1, q, \dots, q^{n-1})$, 
it follows from Lemma~\ref{lem:root} that $d$ does not divide $n$.

In order to show the inclusion $D \subset C$, let $z \in D$ be a primitive root of $1$ 
with $d \mid k$ and $d \nmid n$.
Let $\lambda$ be a partition of $n$.
Since $\sum_i \lambda_i = n$ and $d$ does not divide $n$, there exists a part $\lambda_i$ 
not divisible by $d$.
Then, by using Lemma~\ref{lem:root}, 
$z$ is a root of $h_{\lambda_i} (1, q, \dots, q^{k-1})$, 
so $z$ is a root of $h_\lambda (1, q, \dots, q^{k-1})$.
\end{demo}

\begin{demo}{Proof of Claim~2}
Suppose that $z$ is a common root of $h_\lambda (1, q, \dots, q^{k-1})$'s with $\lambda \vdash n$.
It follows from Claim~1 that $z$ is a primitive $d$-th root of unity for some $d$ 
satisfying $d \mid k$ and $d \nmid n$.
We write $n = s d + t$ with $s \in \Int$ and $0 \le t < d$, 
and consider a partition
$$
\mu = (\underbrace{d, \dots, d}_s, t) \vdash n.
$$
Then
$$
h_\mu (1, q, \dots, q^{k-1})
 = 
h_d (1, q, \dots, q^{k-1})^s \cdot h_t (1, q, \dots, q^{k-1})
$$
and it follows from Lemma~\ref{lem:root} that $z$ is 
a simple root of $h_\mu (1, q, \dots, q^{k-1})$.
\end{demo}

Thus we are done for the proof of (2).

Next we prove (1).
Let $g$ be the greatest common divisor of $s_\lambda(1^k)$ with $\lambda \vdash n$.
Since $\{ e_\lambda : \lambda \vdash n \}$ is another $\Int$-basis of $\Lambda_\Int$, 
the integer $g$ is the greatest common divisor of $e_\lambda(1^k)$'s.

From (2), we already know that $[k]_q/[d]_q$ divides $e_\lambda(1, q, \dots, q^{k-1})$ 
for any partition $\lambda \vdash n$.
Since $[k]_q/[d]_q$ is a monic polynomial, we see that $k/d$ divides $e_\lambda(1^k)$.
Hence $k/d$ divides $g$.
We shall show that $g$ divides $k/d$.

We fix a prime $p$ and show that $\nu_p(g) \le \nu_p(k/d)$.
We put $G = \nu_p(g)$ and $K = \nu_p(k)$.
We write $n = s \cdot p^K + r$ with $s \in \Nat$ and $0 \le r < p^K$, and consider a partition
$$
\lambda = (\underbrace{p^K, \dots, p^K}_s, r) \vdash n.
$$
Then we have
$$
p^G \mid e_\lambda(1^k),
\quad\text{i.e.,}\quad
\left.
p^G \;\middle|\; \binom{k}{p^K}^s \binom{k}{r}
\right.
.
$$
It follows from Lemma~\ref{lem:Kummer} that $\binom{k}{p^K}$ is not divisibly by $p$, so we have
$$
\left.
p^G \;\middle|\; \binom{k}{r}
\right..
$$

If $r = 0$, then $G = 0$ and $\nu_p(g) = 0 \le \nu_p(k/d)$ as desired.
So we may assume $r > 0$.
We put $R = \nu_p(r)$. 
Then it follows from Lemma~\ref{lem:Kummer} that $\nu_p(\binom{k}{r}) = K-R$.
Hence we conclude that $K - R \ge G$.
On the other hand, since $\nu_p(n) = \nu_p(r) = R < K = \nu_p(k)$, we have
$\nu_p( d ) = R$ and $\nu_p(k/d) = K - R$.
This completes the proof of Theorem~\ref{thm:gcd}.
\end{demo}

Let $k$ and $n$ be positive integers and $\lambda$ be a partition of $n$.
Then Theorem~\ref{thm:gcd} (2) implies that 
$s_\lambda(1, q, \dots, q^{k-1})$ is divisible by $[k]_q/[d]_q$, where $d = \gcd(n,k)$.
Since $[k]_q/[d]_q$ is a monic polynomial with integer coefficients, 
the quotient
$$
\frac{ s_\lambda(1, q, \dots, q^{k-1}) }
     { [k]_q/[d]_q }
$$
is also a polynomial with integer coefficients.
In fact, we can show that this quotient has nonnegative integer coefficients.
Recall that a polynomial $f(q) = \sum_{i=0}^m a_i q^i$ of degree $m$ 
is symmetric (resp. unimodal) if $a_i = a_{m-i}$ for all $i$ 
(resp. if there is an index $p$ such that $a_0 \le a_1 \le \dots \le a_{p-1} \le 
a_p \ge a_{p+1} \ge \dots \ge a_d$).
The idea of the following lemma goes back to \cite[Theorem~2]{Andrews}.

\begin{lemma}
\label{lem:Nq}
(See Guo--Krattenthaler \cite[Lemma~5.1]{GK} for example.)
Let $g(q)$ be a polynomial with nonnegative integer coefficients, 
and let $k$ and $d$ be positive integers.
Assume that $g(q)$ is symmetric and unimodal 
and that $[d]_q f(q)/[k]_q$ is a polynomial in $q$.
Then $[d]_q f(q)/[k]_q$ has nonnegative integer coefficients.
\end{lemma}

\begin{prop}
\label{prop:skd_Nq}
Let $k$ and $n$ be positive integers and $d = \gcd(n,k)$.
For a partition $\lambda$ of $n$, the quotient
$$
f_\lambda(q)
 =
\frac{ s_\lambda(1, q, \dots, q^{k-1}) }
     { [k]_q/[d]_q }
$$
is a polynomial with nonnegative integer coefficients.
\end{prop}

This proposition is also given by Garsia--Leven--Wallach--Xin \cite[Theorem~2.1]{GLWX}.

\begin{demo}{Proof}
It is known (\cite[I.8, Example~4]{Macdonald}) 
that $s_\lambda(1, q, \dots, q^{k-1})$ is a symmetric unimodal polynomial 
with nonnegative integer coefficients.
By Theorem~\ref{thm:gcd}, $f_\lambda(q)$ is a polynomial in $q$.
Hence we can apply Lemma~\ref{lem:Nq} to obtain the conclusion.
\end{demo}

\subsection{%
Proof of Theorem~\ref{thm:main} for $W = \Sym_n$
}

Now we are ready to prove Theorem~\ref{thm:main} for $W = \Sym_n$.
Since the implication (i) $\implies$ (ii) follows from Proposition~\ref{prop:Catalan=mult},
it remains to show the implications (ii) $\implies$ (iii), (iii) $\implies$ (i), 
(iv) $\implies$ (iii), and (iii) $\implies$ (v).

\begin{demo}{Proof of (ii) $\implies$ (iii)}
We prove that, if $\Cat_k(\Sym_n, q)$ is a polynomial in $q$, then $\gcd(n,k) = 1$.
We note that
$$
\Cat_k(\Sym_n,q)
 = 
\frac{ h_n(1, q, \dots, q^{k-1}) }
     { [k]_q }
 =
\frac{ h_n(1, q, \dots, q^{k-1}) }
     { \prod_{e \mid k, \, e \neq 1} \Phi_e(q) }.
$$
Let $e (\neq 1)$ be a divisor of $k$.
Since $\Cat_k(\Sym_n,q)$ is a polynomial, 
any $e$-th primitive root of $1$ must be a root of $h_n(1, q, \dots, q^{k-1})$.
It follows from Lemma~\ref{lem:root} that $e$ does not divide $n$.
Hence we see that $k$ is coprime to $n$.
\end{demo}

\begin{demo}{Proof of (iii) $\implies$ (i)}
Suppose that $\gcd(n,k) = 1$.
By (\ref{eq:decomp_S1}), it is enough to show that $s_\lambda(1, q, \dots, q^{k-1})/[k]_q$ 
is a polynomial with nonnegative integer coefficients for all $\lambda \vdash n$.
This is a consequence of Proposition~\ref{prop:skd_Nq}.
\end{demo}

\begin{demo}{Proof of (iii) $\implies$ (v)}
If $\gcd(n,k) = 1$, then $m_\lambda(1^k)/k \in \Nat$ by Theorem~\ref{thm:gcd} (2).
Hence it follows from (\ref{eq:decomp_S3}) that $\varphi_k$ is a permutation character.
\end{demo}

\begin{demo}{Proof of (iv) $\implies$ (iii)}
If $\varphi_k$ is a character of some representation, 
then it follows from (\ref{eq:decomp_S2}) 
that $s_\lambda(1^k)/k$ is an integer for all $\lambda \vdash n$.
Then $k$ divides $\gcd \{ s_\lambda(1^k) : \lambda \vdash n \}$.
Hence, by using Theorem~\ref{thm:gcd} (1), we see that $k$ divides $k/\gcd(n,k)$.
So we must have $\gcd(n,k) = 1$.
\end{demo}

This completes the proof of Theorem~\ref{thm:main} for the symmetric groups.
\section{%
$G(m,p,n)$
}

This section is devoted to the proof of the main result for the groups of the form $W = G(m,p,n)$, 
except for the dihedral groups.
In this section, we fix an integer $m \ge 2$ and put $\zeta = e^{2 \pi \sqrt{-1}/m}$.

\subsection{%
Representations of $G(m,1,n)$ and symmetric functions
}

We review several facts concerning representations of 
$G(m,1,n)$ and their connection to symmetric functions.
See \cite[Chapter~1, Appendix B]{Macdonald} for details.

Let $\Par^{(m)}_n$ be the set of $m$-tuples of partitions 
$\boldsymbol{\lambda} = (\lambda^{(0)}, \cdots, \lambda^{(m-1)})$ with 
$| \boldsymbol{\lambda} | = \sum_{i=0}^{m-1} |\lambda^{(i)}| = n$.
We put $\Par^{(m)} = \bigsqcup_{n \ge 0} \Par^{(m)}_n$.

The group $G(m,1,n)$ is realized as a complex reflection group acting on $V = \Comp^n$, 
which consists of all monomial matrices such that the nonzero entries are $m$-th roots of unity.
In this subsection we write $\Sym^{(m)}_n = G(m,1,n)$ for short.
The group $\Sym^{(m)}_n$ is isomorphic to the semidirect product 
$\Sym_n \ltimes C_m^n$, where $C_m$ is the cyclic group of order $m$.
 
It is known that an element $w \in \Sym^{(m)}_n$ is conjugate to the matrix of the form
$$
\bigoplus_{i=0}^{m-1} \bigoplus_{j \ge 1} C(\mu^{(i)}_j,\zeta^i),
$$
where $\boldsymbol{\mu} = (\mu^{(0)}, \dots, \mu^{(m-1)}) \in \Par^{(m)}_n$, 
called the cycle type of $w$ and denoted by $\type(w)$, 
and $C(l,\alpha)$ is the $l \times l$ matrix given by
$$
C(l,\alpha)
 =
\begin{pmatrix}
 0 & 0 & 0 & \cdots & 0 & \alpha \\
 1 & 0 & 0 & \cdots & 0 & 0 \\
 0 & 1 & 0 & \cdots & 0 & 0 \\
 \vdots & \vdots & \vdots & \ddots & \vdots & \vdots \\
 0 & 0 & 0 & \cdots & 0 & 0 \\
 0 & 0 & 0 & \cdots & 1 & 0
\end{pmatrix}.
$$
And two elements of $\Sym^{(m)}_n$ are conjugate if and only if they have the same cycle types.
Hence the conjugacy classes of $\Sym^{(m)}_n$ are indexed by $\Par^{(m)}_n$.
It is easy to see that, if $w \in \Sym^{(m)}_n$ have the cycle type $\boldsymbol{\mu}$, then
\begin{equation}
\label{eq:det_G}
\det{}_V (1 - t w)
 =
\prod_{i=0}^{m-1} \prod_{j \ge 0} \left( 1 - \zeta^i t^{\mu^{(i)}_j} \right).
\end{equation}

The irreducible characters of $\Sym^{(m)}_n$ are also parametrized by $\Par^{(m)}_n$.
We denote by $\chi^{\boldsymbol{\lambda}}$ the irreducible character of $\Sym^{(m)}_n$ 
corresponding to $\boldsymbol{\lambda}$, 
and by $\chi^{\boldsymbol{\lambda}}_{\boldsymbol{\mu}}$ the character value 
of $\chi^{\boldsymbol{\lambda}}$ at an element of cycle type $\boldsymbol{\mu}$.

Let $R(\Sym^{(m)}_n)$ be the vector space of complex-valued class functions on $\Sym^{(m)}_n$, 
and put $R(\Sym^{(m)}_\bullet) = \bigoplus_{n \ge 0} R(\Sym^{(m)}_n)$.
Then $R(\Sym^{(m)}_\bullet)$ becomes a graded $\Comp$-algebra with respect to the product
$$
f \cdot g = \Ind_{\Sym^{(m)}_n \times \Sym^{(m)}_l}^{\Sym^{(m)}_{n+l}} (f \times g),
$$
where $f \in R(\Sym^{(m)}_n)$ and $g \in R(\Sym^{(m)}_l)$.

Let $\Lambda(\vectx^{(i)})$ be the ring of symmetric functions in variables $\vectx^{(i)}
 = (x^{(i)}_1, x^{(i)}_2, \dots)$, and put
$$
\Lambda^{(m)}
 =
\Lambda(\vectx^{(0)}) \otimes \Lambda(\vectx^{(1)}) \otimes \cdots \otimes \Lambda(\vectx^{(m-1)}).
$$
For $\boldsymbol{\lambda} \in \Par^{(m)}_n$, we put
\begin{align*}
S_{\boldsymbol{\lambda}}
 &=
s_{\lambda^{(0)}}(\vectx^{(0)}) 
s_{\lambda^{(1)}}(\vectx^{(1)}) 
\dots
s_{\lambda^{(m-1)}}(\vectx^{(m-1)}),
\\
P_{\boldsymbol{\lambda}}
 &=
\prod_{i=0}^{m-1} \prod_{j \ge 1}
 \left(
  \sum_{r=0}^{m-1} \zeta^{ir} p_{\mu^{(i)}_j}(\vectx^{(r)})
 \right).
\end{align*}
We define the Frobenius characteristic map $\ch^{(m)} : R(\Sym^{(m)}_\bullet) \to \Lambda^{(m)}$ by
$$
\ch^{(m)}( f )
 =
\frac{ 1 }{ m^n n! }
\sum_{w \in \Sym^{(m)}_n} 
 f(w) 
 \overline{ P_{\type(w)} }
\quad(f \in R(\Sym^{(m)}_n),
$$
where
$$
\overline{ P_{\boldsymbol{\lambda}} }
 = 
\prod_{i=0}^{m-1} \prod_{j \ge 1}
 \left(
   \sum_{r=0}^{m-1} \zeta^{-ir} p_{\mu^{(i)}_j}(\vectx^{(r)})
 \right).
$$
Then we have

\begin{prop}
\label{prop:Frobenius_G}
\begin{enumerate}
\item[(1)]
The Frobenius characteristic map $\ch^{(m)}$ is an algebra isomorphism.
\item[(2)]
For $\boldsymbol{\lambda} \in \Par^{(m)}_n$, we have
\begin{equation}
\label{eq:Frobenius_G1}
\ch^{(m)} ( \chi^{\boldsymbol{\lambda}} )
 =
S_{\boldsymbol{\lambda}}.
\end{equation}
\item[(3)]
For $\boldsymbol{\mu} \in \Par^{(m)}_n$, we have
\begin{equation}
\label{eq:Frobenius_G2}
P_{\boldsymbol{\mu}}
 =
\sum_{\boldsymbol{\lambda} \in \Par^{(m)}_n}
 \chi^{\boldsymbol{\lambda}}_{\boldsymbol{\mu}} 
 S_{\boldsymbol{\lambda}}.
\end{equation}
\end{enumerate}
\end{prop}

We can use this setting to give a formula for the multiplicities of irreducible characters 
in the class function $\hat{\varphi}$ given by (\ref{def:phiqu}).

\begin{prop}
\label{prop:decomp_G1}
For $\boldsymbol{\lambda} \in \Par^{(m)}_n$, 
the multiplicity $\hat{m}^{\boldsymbol{\lambda}}$ of $\chi^{\boldsymbol{\lambda}}$ 
in $\hat{\varphi}$ is given by
$$
\hat{m}^{\boldsymbol{\lambda}}
 =
\hat{\pi} \left( S_{\boldsymbol{\lambda}} \right),
$$
where $\hat{\pi} : \Lambda^{\otimes m} \to \Comp(q,u)$ is the ring homomorphism given by
$$
\hat{\pi} ( p_l(\vectx^{(r)}) ) 
=
\begin{cases}
\dfrac{ 1 - q^{(m-1)l} u^l }{ 1 - q^{ml} } &\text{if $r=0$,} \\
q^{rl} \dfrac{ 1 - q^{-l} u^l }{ 1 - q^{ml} } &\text{if $1 \le r \le m-1$.}
\end{cases}
$$
\end{prop}

\begin{demo}{Proof}
By a direct computation, we have
$$
\frac{ 1 - q^{(m-1)l} u^l }{ 1 - q^{ml} }
 +
\sum_{s=1}^{m-1}
 \zeta^{is} q^{sl} \frac{ 1 - q^{-l} u^l }{ 1 - q^{ml} }
 =
\frac{ 1 - \zeta^i u^l }{ 1 - \zeta^i q^l }.
$$
Hence we see that
$$
\hat{\varphi}(w)
 =
\hat{\pi} \left( P_{\type(\mu)} \right)
\quad(w \in G(m,1,n)).
$$
It follows from the Frobenius formula (\ref{eq:Frobenius_G2}) that
$$
\hat{\varphi}(w)
 =
\sum_{\lambda \in \Par^{(m)}_n} 
 \hat{\pi} \left( S_{\boldsymbol{\lambda}} \right)
 \chi^{\boldsymbol{\lambda}}(w)
\quad(w \in G(m,1,n)).
$$
\end{demo}

\begin{corollary}
\label{cor:decomp_G2}
\begin{enumerate}
\item[(1)]
The multiplicity $m_k^{\boldsymbol{\lambda}}$ of $\chi^{\boldsymbol{\lambda}}$ 
in $\varphi_k$ is given by
\begin{equation}
\label{eq:decomp_G1}
m_k^{\boldsymbol{\lambda}}
 =
\prod_{x \in \lambda^{(0)}}
 \frac{ \frac{m+k-1}{m} + c(x) }
      { h(x) }
\prod_{r=1}^{m-1}
 \prod_{x \in \lambda^{(r)}}
 \frac{ \frac{k-1}{m} + c(x) }
      { h(x) },
\end{equation}
where $x$ runs over all cells in the Young diagram of $\lambda^{(0)}$ 
or $\lambda^{(r)}$, and $h(x)$ and $c(x)$ denote the hook length and the content of $x$ 
respectively.
\item[(2)]
If $k = pm + 1$ with $p \in \Int$, 
then the multiplicity $\tilde{m}^{\boldsymbol{\lambda}}$ of $\chi^{\boldsymbol{\lambda}}$ 
in $\tilde{\varphi}_k$ is given by
\begin{equation}
\label{eq:decomp_G2}
\tilde{m}^{\boldsymbol{\lambda}}_k
 =
s_{\lambda^{(0)}} (1, q^m, \cdots, q^{pm} )
\prod_{r = 1}^{m-1}
 s_{\lambda^{(r)}} (q^r, q^{r+m}, \cdots, q^{r+(p-1)m} ),
\end{equation}
which is a polynomial in $q$ with nonnegative integer coefficients.
\end{enumerate}
\end{corollary}

\begin{demo}{Proof}
By replacing $u$ by $q^k$ in Proposition~\ref{prop:decomp_G1}, we see that 
the multiplicity $\tilde{m}^{\boldsymbol{\lambda}}_k$ in $\tilde{\varphi}_k$ is given by
$$
\tilde{m}^{\boldsymbol{\lambda}}_k = \tilde{\pi}_k \left( S_{\boldsymbol{\lambda}} \right),
$$
where the algebra homomorphism $\tilde{\pi}_k : \Lambda^{(m)} \to \Comp(q)$ is defined by 
$$
\tilde{\pi}_k ( p_l(\vectx^{(r)}) ) 
=
\begin{cases}
\dfrac{ 1 - q^{(m+k-1)l} }{ 1 - q^{ml} } &\text{if $r=0$,} \\
q^{rl} \dfrac{ 1 - q^{(k-1)l} }{ 1 - q^{ml} } &\text{if $1 \le r \le m-1$.}
\end{cases}
$$

(1)
By putting $q=1$, we see that the the multiplicity $m^{\boldsymbol{\lambda}}_k$ in $\varphi_k$ is given by 
the specialization
$$
m^{\boldsymbol{\lambda}}_k = \pi_k \left( S_{\boldsymbol{\lambda}} \right),
$$
where the algebra homomorphism $\pi_k : \Lambda^{(m)} \to \Comp(q)$ is defined by 
$$
\pi_k ( p_l(\vectx^{(r)}) ) 
=
\begin{cases}
\dfrac{ m+k-1 }{ m } &\text{if $r=0$,} \\
\dfrac{ k-1 }{ m } &\text{if $1 \le r \le m-1$.}
\end{cases}
$$
Hence (\ref{eq:decomp_G1}) follows from the specialization formula 
(\cite[I.3 Example~4]{Macdonald})
$$
\rho_z (s_\lambda)
 = 
\prod_{x \in \lambda}
 \frac{ z + c(x) }
      { h(x) },
$$
where $\pi_z : \Lambda \to \Comp(z)$ is defined by $\pi_z(p_l) = z$ for $l \ge 1$.

(2)
If $k = pm+1$, then the specialization $\tilde{\pi}_k$ is given by substitution
$$
x^{(r)}_j
 =
\begin{cases}
 q^{(j-1)m} &\text{if $r=0$ and $1 \le j \le p+1$,} \\
 q^{r+(j-1)m} &\text{if $1 \le r \le m-1$ and $1 \le j \le p$,} \\
 0 &\text{otherwise.}
\end{cases}
$$
Hence we obtain (\ref{eq:decomp_G2}).
\end{demo}

Next we give an expression of $\varphi_k$ as a linear combination of permutation characters.

\begin{lemma}
\label{lem:StoG}
Let $F : R(\Sym_\bullet) \to R(\Sym^{(m)}_\bullet)$ be the linear map defined by
$$
F(f) = \Ind_{\Sym_n}^{\Sym^{(m)}_n} f
\quad(f \in R(\Sym_n)).
$$
Then $F$ is an graded algebra homomorphism and
the following diagram commutes.
$$
\begin{CD}
R(\Sym_\bullet) @>{F}>>      R(\Sym^{(m)}_\bullet) \\
@V{\ch}VV                    @VV{\ch^{(m)}}V \\
\Lambda         @>>{\Delta}> \Lambda^{(m)}
\end{CD}
$$
where $\Delta : \Lambda \to \Lambda^{\otimes m}$ is the coproduct given by the plethystic substitution
$$
\Delta(f) = f[\vectx^{(0)} + \vectx^{(1)} + \cdots + \vectx^{(m-1)}].
$$
\end{lemma}

\begin{demo}{Proof}
By using the definition of the products on $R(\Sym_\bullet)$ and $R(\Sym^{(m)}_\bullet)$, 
we see that
$$
F(f \cdot g)
 = \Ind_{\Sym_n \times \Sym_l}^{\Sym^{(m)}_{n+l}} (f \times g)
 = F(f) \cdot F(g)
$$
for $f \in R(\Sym_n)$ and $g \in R(\Sym_l)$.

If we denote by $c_\mu$ the characteristic function of the conjugacy class of $\Sym_n$ 
corresponding to a partition $\mu$, then we have
$$
\ch c_\mu = \frac{1}{z_\mu} p_\mu,
$$
where $z_\mu = \prod_{i \ge 1} i^{m_i} m_i!$.
And $F(c_\mu)$ is the characteristic function of the conjugacy class of $\Sym^{(m)}_n$ 
corresponding to $(\mu, \emptyset, \dots, \emptyset)$, so we have
$$
\ch^{(m)} (F(c_\mu))
 = 
\frac{1}{z_\mu} \prod_{i=1}^{l(\mu)}
 \left( p_{\mu_i}(\vectx^{(0)}) + \dots + p_{\mu_i}(\vectx^{(m-1)}) \right)
 =
\frac{1}{z_\mu} p_\mu[ \vectx^{(0)} + \dots + \vectx^{(m-1)}].
$$
Hence the diagram commutes.
\end{demo}

\begin{prop}
\label{prop:decomp_G3}
If $k = pm+1$ with $p \in \Int$, then we have
\begin{equation}
\label{eq:decomp_G3}
\varphi^{\Sym^{(m)}_n}_k
 =
\sum_{r=0}^n \sum_{\lambda \vdash n-r}
 m_\lambda(1^p) \eta^{r,\lambda}.
\end{equation}
where $\eta^{r,\lambda}$ is the permutation character on 
the coset space $\Sym^{(m)}_n / \Sym^{(m)}_r \times \Sym_\lambda$ 
by the parabolic subgroup $\Sym^{(m)}_r \times \Sym_\lambda
 = \Sym^{(m)}_r \times \Sym_{\lambda_1} \times \Sym_{\lambda_2} \times \dots$.
\end{prop}

\begin{demo}{Proof}
If $k = pm+1$, then it follows from (\ref{eq:decomp_G2}) with $q=1$ that
$$
\ch^{(m)} ( \varphi^{\Sym^{(m)}_n}_k )
=
\sum_{ \boldsymbol{\lambda} \in \Par^{(m)}_n }
 s_{\lambda^{(0)}} (1^{p+1}) \prod_{i=1}^{m-1} s_{\lambda^{(i)}} (1^p)
 \prod_{i=0}^{m-1} s_{\lambda^{(i)}} (\vectx^{(i)}).
$$
Hence, by using the Cauchy identity (\ref{eq:Cauchy}), we have
\begin{align*}
\sum_{n \ge 0} \ch ( \varphi^{\Sym^{(m)}_n}_k )
 &=
\sum_{ \boldsymbol{\lambda} \in \Par^{(m)} }
 s_{\lambda^{(0)}} (1^{p+1}) \prod_{i=1}^{m-1} s_{\lambda^{(i)}} (1^p)
 \prod_{i=0}^{m-1} s_{\lambda^{(i)}} (\vectx^{(i)})
\\
 &=
\prod_j (1 - x^{(0)}_j)^{-p-1}
\prod_{i=1}^{m-1} \prod_j (1 - x^{(i)}_j)^{-p}
\\
 &=
\prod_j (1 - x^{(0)}_j)^{-1}
\prod_{i=0}^{m-1} \prod_j (1 - x^{(i)}_j)^{-p}
\\
 &=
\left( \sum_{r=0}^\infty h_r(\vectx^{(0)}) \right)
\left(
  \sum_\lambda m_\lambda(1^p) 
  h_\lambda \left[ \vectx^{(0)} + \vectx^{(1)} + \cdots + \vectx^{(r-1)} \right]
\right).
\end{align*}
Therefore we have
$$
\ch ( \varphi^{\Sym^{(m)}_n}_k )
 =
\sum_{r=0}^n \sum_{\lambda \vdash n-r} 
 m_\lambda(1^p) h_r(\vectx^{(0)}) \Delta(h_\lambda).
$$
Here $h_r(\vectx^{(0)}) = S_{((r), \emptyset, \dots, \emptyset)}$ is 
the Frobenius characteristic of the trivial representation of $\Sym^{(m)}_r$.
And it follows from Lemma~\ref{lem:StoG} that $\Delta(h_\lambda)$ is the Frobenius characteristic 
of the permutation character on $\Sym^{(m)}_{n-r}/\Sym_\lambda$.
Hence we have
$$
\varphi^{\Sym^{(m)}_n}_k
 =
\sum_{r=0}^n \sum_{\lambda \vdash n-r}
 m_\lambda(1^p) 
 \Ind_{\Sym^{(m)}_r \times \Sym_\lambda}^{\Sym^{(m)}_n} \triv.
$$
\end{demo}

\subsection{%
Representation Theory of $G(m,p,n)$
}

For a divisor $p$ of $m$, the group $G(m,p,n)$ is the subgroup of $G(m,1,n)$ 
consisting of the matrices in $G(m,1,n)$ such that the product of the nonzero entries 
is $(m/p)$-th root of unity.
Then $G(m,p,n)$ is a normal subgroup of $G(m,1,n)$ of index $p$, 
and the quotient group $G(m,1,n)/G(m,p,n)$ is the cyclic group of order $p$.
Hence we can apply Clifford Theory to obtain irreducible representations of $G(m,p,n)$ 
from those of $G(m,1,n)$.
See \cite[Section~6]{Stembridge} for details.

Let $\sh : \Par^{(m)}_n \to \Par^{(m)}_n$ be the shift operator defined by
$$
\sh( \lambda^{(0)}, \cdots, \lambda^{(m-1)} )
 =
( \lambda^{(1)}, \cdots, \lambda^{(m-1)}, \lambda^{(0)} ).
$$
Then the group $\langle \sh^{m/p} \rangle$ acts on $\Par^{(m)}_n$.
For $\boldsymbol{\lambda} \in \Par^{(m)}_n$, we have
$$
\Res^{G(m,1,n)}_{G(m,p,n)}
\chi^{\boldsymbol{\lambda}}
=
\Res^{G(m,1,n)}_{G(m,p,n)}
\chi^{\sh^{m/p}(\boldsymbol{\lambda})}.
$$
If the stabilizer of $\chi^{\boldsymbol{\lambda}}$ in $\langle \sh^{m/p} \rangle$ 
has the order $t$, then $\Res^{G(m,1,n)}_{G(m,p,n)} \chi^{\boldsymbol{\lambda}}$ 
is decomposed into the sum of $t$ distinct irreducible characters with multiplicity $1$.
And any irreducible character of $G(m,p,n)$ can be obtained in this way.
Since $\varphi^{G(m,p,n)}_k$ is the restriction of $\varphi^{G(m,1,n)}_k$ to 
$G(m,p,n)$, we obtain the following proposition.

\begin{prop}
\label{prop:mult_Gp}
If an irreducible character $\chi$ of $G(m,p,n)$ appears in the restriction 
$\Res^{G(m,1,n)}_{G(m,p,n)} \chi^{\boldsymbol{\lambda}}$ with 
$\boldsymbol{\lambda} \in \Par^{(m)}_n$, 
then the multiplicity of $\chi$ in $\varphi_k^{G(m,p,n)}$ is equal to
$$
\sum_{ \boldsymbol{\mu} \in \langle \sh^{m/p} \rangle \boldsymbol{\lambda} }
 m^{\boldsymbol{\mu}}_k
$$
where the summation is taken over the $\langle \sh^{m/p} \rangle$-orbit of 
$\boldsymbol{\lambda}$ and 
$m^{\boldsymbol{\mu}}_k$ is the multiplicity of $\chi^{\boldsymbol{\mu}}$ 
in $\varphi^{G(m,1,n)}_k$.
\end{prop}

\subsection{%
Proof of Theorem~\ref{thm:main} for $G(m,p,n)$ with $n \ge 3$
}

In this subsection, we assume $n \ge 3$.

\begin{demo}{Proof of (ii) $\implies$ (iii)}
We show that, if $\Cat_k(W,q)$ is a polynomial, then $k \equiv 1 \bmod m$.
Since the degrees $(d_1, \dots, d_n)$ of $W = G(m,n,p)$ are given by
$$
(d_1, \dots, d_n)
 =
(m, 2m, \dots, (n-1)m, mn/p),
$$
we have
$$
\Cat_k(W, q)
 =
\frac{ \prod_{i=1}^{n-1} [k+im-1]_q \cdot [k+mn/p-1]_q }
     { \prod_{i=1}^{n-1} [im]_q \cdot [mn/p]_q }.
$$
Since $n \ge 3$, the denominator of $\Cat_k(W,q)$ is divisible by $\Phi_m(q)^2$, 
where $\Phi_m(q)$ is the $m$-th cyclotomic polynomial.
So, if $\Cat_k(W,q)$ is a polynomial, then the numerator is also divisible by $\Phi_m(q)^2$.
Hence at least one of the factors $[k+im-1]_q$ with $1 \le i \le n-1$ is divisible by $\Phi_m(q)$.
This implies $m \mid k+im-1$, so $k \equiv 1 \bmod m$.
\end{demo}

\begin{demo}{Proof of (iii) $\implies$ (i)}
Suppose that $k \equiv 1 \bmod m$.
Then it follows from (\ref{eq:decomp_G2}) that 
$\tilde{\varphi}^{G(m,1,n)}_k$ is the graded character of a graded representation.
By restricting to $G(m,p,n)$, we see that $\tilde{\varphi}^{G(m,p,n)}_k$ is also 
the graded character of a graded representation.
\end{demo}

\begin{demo}{Proof of (iii) $\implies$ (v)}
Suppose that $k \equiv 1 \bmod m$.
Then it follows from (\ref{eq:decomp_G3}) that
$\varphi^{G(m,1,n)}_k$ is a permutation character of $G(m,1,n)$..

In general, by using Mackey's Restriction Theorem (see \cite[Proposition~22]{Serre})
$$
\Res^G_K \Ind_H^G \triv
 =
\sum_g \Ind_{gHg^{-1} \cap K}^K \triv,
$$
where $g$ runs over all double coset representatives of $K \backslash G / H$,
we can see that $\varphi^{G(m,p,n)}_k = \Res^{G(m,1,n)}_{G(m,p,n)} \varphi^{G(m,1,n)}_k$ 
is also a permutation character.
\end{demo}

The most subtle part of the proof is the implication (iv) $\implies$ (iii).

\begin{demo}{Proof of (iv) $\implies$ (iii)}
Suppose that $\varphi^{G(m,p,n)}_k$ comes from a genuine representation of $G(m,p,n)$.
Then, by restricting to $G(m,m,n)$, we see that $\varphi^{G(m,m,n)}_k$ also comes from 
a genuine representation.

Since the degrees and the codegrees of $G(m,m,n)$ are given by
\begin{align*}
(d_1, \dots, d_n)
 &= 
(m, 2m, \dots, (n-1)m, n),
\\
(d^*_1, \dots, d^*_n)
 &=
(0, m, \dots, (n-2)m, (n-1)m-n),
\end{align*}
we can use Proposition~\ref{prop:Catalan=mult} to see that
the multiplicities of $\triv$ and $\det$ in $\tilde{\varphi}^{G(m,m,n)}_k$ are given by
\begin{align*}
m^{\triv}_k
 &=
\frac{1}{m^{n-1} n!} \prod_{i=1}^{n-1} (k + im -1) \cdot (k+n-1),
\\
m^{\det}_k
 &=
\frac{1}{m^{n-1} n!} \prod_{i=1}^{n-1} (k - im + m -1) \cdot (k-(n-1)m+(n-1)).
\end{align*}
This can be also derived from (\ref{eq:decomp_G1}) and Proposition~\ref{prop:mult_Gp}.
We prove that, if $m^{\triv}_k$ and $m^{\det}_k$ are both integers, then $k \equiv 1 \bmod m$.

We consider two polynomials
\begin{align*}
f(z)
 &=
\frac{1}{n!} \prod_{i=1}^{n-1} (z+i) \cdot (mz+n),
\\
g(z)
 &=
\frac{1}{n!} \prod_{i=1}^{n-1} (z-i+1) \cdot (m(z-n+1)+n).
\end{align*}
Then we have
$$
m^{\triv}_k = f \left( \frac{k-1}{m} \right),
\quad
m^{\det}_k = g \left( \frac{k-1}{m} \right).
$$
Hence $z = (k-1)/m$ is a solution of the polynomial equation
\begin{multline*}
\frac{n!}{2} ( f(z) + g(z) ) - (n-2)! z ( f(z) - g(z) )
\\
= \frac{n!}{2} (m^{\triv}_k + m^{\det}_k) - (n-2)! (m^{\triv}_k - m^{\det}_k) z.
\end{multline*}
Recall that a rational number is a rational integer if and only if it is an algebraic number.
Hence it is enough to show that
$$
h(z) = \frac{n!}{2} ( f(z) + g(z) ) - (n-2)! z ( f(z) - g(z) )
$$
is a monic polynomial of degree $n-1$ with integer coefficients.

The polynomials $f(z)$ and $g(z)$ are written as
$$
f(z) = f_1(z) + f_2(z),
\quad
g(z) = g_1(z) + g_2(z),
$$
where
\begin{alignat*}{2}
f_1(z)
 &=
\frac{m}{n!} \prod_{i=1}^n (z+i-1),
&\quad
f_2(z)
 &=
\frac{1}{(n-1)!} \prod_{i=1}^{n-1} (z+i),
\\
g_1(z)
 &=
\frac{m}{n!} \prod_{i=1}^n (z-i+1)
&\quad
g_2(z)
 &=
\frac{1}{(n-1)!} \prod_{i=1}^{n-1} (z-i+1).
\end{alignat*}
Since we have
\begin{align*}
f_1(z)
&=
\frac{m}{n!} \left( z^n + \frac{n(n-1)}{2} z^{n-1} + \text{lower terms} \right),
\\
g_1(z)
&=
\frac{m}{n!} \left( z^n - \frac{n(n-1)}{2} z^{n-1} + \text{lower terms} \right),
\\
f_2(z)
 &=
\frac{1}{(n-1)!} \left( z^{n-1} + \frac{n(n-1)}{2} z^{n-2} + \text{lower terms} \right),
\\
g_2(z)
 &=
\frac{1}{(n-1)!} \left( z^{n-1} - \frac{(n-1)(n-2)}{2} z^{n-2} + \text{lower terms} \right),
\end{align*}
we see that $h(z)$ is a monic polynomial of degree $n-1$.

Next we show that $h(z)$ has integer coefficients.
Let $c(n,j)$ be the signless Stirling number of the first kind.
Then we have (see \cite[1.3.7 Proposition]{EC1} for example)
$$
\prod_{i=1}^n (z+i-1) = \sum_{k=0}^n c(n,j) z^j.
$$
By using this generating function, we see that
\begin{gather*}
\frac{n!}{2} (f_1(u) + g_1(u))
=
m \sum_{\text{$n-j$ is even}} c(n,j) u^j,
\\
(n-2)! (f_1(u) - g_1(u))
=
\frac{ m }{ \binom{n}{2} } \sum_{\text{$n-j$ is odd}} c(n,j) u^j.
\end{gather*}
Also by using the recurrence relation (see \cite[1.3.6 Lemma]{EC1} for example)
$$
c(n,j) = c(n-1,j-1) + (n-1) c(n-1,j),
$$
we have
\begin{align*}
\frac{n!}{2} (f_2(u) + g_2(u))
&=
n \sum_{\text{$n-1-j$ is even}} c(n-1,j) u^j 
+
\binom{n}{2} \sum_{j=1}^n c(n-1,j) u^{j-1},
\\
(n-2)! (f_2(u) - g_2(u))
&=
n \cdot \frac{1}{ \binom{n-1}{2} } \sum_{\text{$n-1-j$ is odd}} c(n-1,j) u^j
+
\sum_{j=1}^n c(n-1,j) u^{j-1}.
\end{align*}
Therefore it is enough to show that $c(n,j)$ is divisible by $\binom{n}{2}$ if $n-j$ is odd.
The proof is reduced to the following lemma.
\end{demo}

\begin{lemma}
\label{lem:Stirling}
(Stanley \cite[Corollary~3.4]{Stanley2011}, Zagier \cite[Application~3]{Zagier})
If $n-j$ is odd, then $c(n,j)$ is divisible by $\binom{n}{2}$.
\end{lemma}

We will give a refinement and another proof in Appendix.

\subsection{%
Proof of Theorem~\ref{thm:main} for $G(m,p,n)$ with $n=2$ and $p < m$
}

In this subsection, we deal with the case where $n=2$ and $p < m$.
The remaining case where $n=2$ and $p=m$ (i.e., $W = G(m,m,2)$ is the dihedral group 
of order $2m$) will be treated in the next section.

The proofs of the implications (iii) $\implies$ (i) and (iii) $\implies$ (v) are 
the same as in the case where $n \ge 3$.

\begin{demo}{Proof of (ii) $\implies$ (iii)}
The degrees and the codegrees of $G(m,p,2)$ are given by
$$
(d_1, d_2) = (m, mn/p),
\quad
(d^*_1, d^*_2) = (0, m),
$$
so we have
$$
\Cat^*_k(W,q)
 = 
q^{m+2}
\frac{ [k-1]_q [k-m-1]_q }
     { [m]_q [mn/p]_q }.
$$
If it is a polynomial in $q$, 
then the numerator is divisible by $\Phi_m(q)$, so
$m \mid (k-1)$ or $m \mid (k-m-1)$.
Hence we have $k \equiv 1 \bmod m$.
\end{demo}

\begin{demo}{Proof of (iv) $\implies$ (iii)}
Let $\chi$ and $\eta$ be the irreducible characters of $G(m,p,2)$ corresponding 
to the $\langle \sh^{m/p} \rangle$-orbits of
$$
\boldsymbol{\lambda} = ((2), \emptyset, \cdots, \emptyset),
\quad\text{and}\quad
\boldsymbol{\mu} = (\emptyset, (2), \emptyset, \cdots, \emptyset),
$$
respectively. Since $p < m$, these two orbits are distinct, so $\chi \neq \eta$.
Then, by using (\ref{eq:decomp_G1}) and Proposition~\ref{prop:mult_Gp}, 
the multiplicities of $\chi$ and $\eta$ in $\varphi_k$ are given by
\begin{align*}
m^{\chi}
 &=
\frac{1}{2 m^2} (k+m-1)(k+2m-1)
+ (p-1) \cdot \frac{1}{2 m^2} (k-1)(k+m-1)
\\
 &=
\frac{p}{2 m^2} (k+m-1) \left( k + \frac{2m}{p} - 1 \right),
\\
m^{\eta}
 &=
\frac{p}{2 m^2} (k-1)(k+m-1).
\end{align*}
So we have
$$
m^{\chi} - m^{\eta}
 =
\frac{k-1}{m}.
$$
Thus, if both of the multiplicities $m^\chi$ and $m^\eta$ are integers, 
then $k \equiv 1 \bmod m$.
\end{demo}
\section{%
Dihedral groups
}

This section is devoted to the proof of Theorem~\ref{thm:main} for the dihedral groups.

\subsection{%
Preliminaries
}

Let $m$ be a positive integer.
The dihedral group $D_{2m}$ of order $2m$ is defined by generators and relations:
$$
D_{2m} = \langle a, b \mid a^m = b^2 = 1, \, b^{-1} a b = a^{-1} \rangle.
$$
The dihedral group $D_{2m}$ is a complex reflection group acting on $\Comp^2$ via
$$
a
 = 
\begin{pmatrix} 
\cos 2 \pi/m & - \sin 2 \pi/m \\
\sin 2 \pi/m & \cos 2 \pi/m
\end{pmatrix},
\quad
b
 =
\begin{pmatrix}
 -1 & 0 \\
 0 & 1
\end{pmatrix}.
$$
Then the conjugacy classes and the irreducible characters of $D_{2m}$ are described as follows.
(See a standard textbook on representation theory of finite groups, 
e.g. \cite[5.3]{Serre}.)

\begin{prop}
\label{prop:character_D}
\begin{enumerate}
\item[(1)]
Suppose $m$ is even.
Then the $(m/2+3)$ elements
$$
1, \quad a^i \quad(1 \le i \le m/2), \quad b, \quad ab
$$
form a complete set of representatives of conjugacy classes of $D_{2m}$.
And $D_{2m}$ has four one-dimensional irreducible characters $\xi_0, \xi_1, \xi_2, \xi_3$ 
and $m/2-1$ two-dimensional irreducible characters $\chi_j$ ($1 \le j \le m/2-1$), 
whose character values are given by
\begin{alignat*}{4}
\xi_0(1) &= 1, \quad&
\xi_0(a^i) &= 1, \quad&
\xi_0(b) &= 1, \quad&
\xi_0(ab) &= 1,
\\
\xi_1(1) &= 1, \quad&
\xi_1(a^i) &= 1, \quad&
\xi_1(b) &= -1, \quad&
\xi_1(ab) &= -1,
\\
\xi_2(1) &= 1, \quad&
\xi_2(a^i) &= (-1)^i, \quad&
\xi_2(b) &= 1, \quad&
\xi_2(ab) &= -1,
\\
\xi_3(1) &= 1, \quad&
\xi_3(a^i) &= (-1)^i, \quad&
\xi_3(b) &= -1, \quad&
\xi_3(ab) &= 1,
\\
\chi_j(1) &= 2, \quad&
\chi_j(a^i) &= \zeta^{ij}+\zeta^{-ij}, \quad&
\chi_j(b) &= 0, \quad&
\chi_j(ab) &= 0.
\end{alignat*}
\item[(2)]
Suppose $m$ is odd.
Then the $((m-1)/2+2)$ elements
$$
1, \quad a^i \quad (1 \le i \le (m-1)/2), \quad b
$$
form a complete set of representatives of conjugacy classes of $D_{2m}$.
And $D_{2m}$ has two one-dimensional irreducible characters $\xi_0, \xi_1$ 
and $(m-1)/2$ two-dimensional irreducible characters $\chi_j$ ($1 \le j \le (m-1)/2$), 
whose character values are given by
\begin{alignat*}{3}
\xi_0(1) &= 1, \quad&
\xi_0(a^i) &= 1, \quad&
\xi_0(b) &= 1,
\\
\xi_1(1) &= 1, \quad&
\xi_1(a^i) &= 1, \quad&
\xi_1(b) &= -1,
\\
\chi_j(1) &= 2, \quad&
\chi_j(a^i) &= \zeta^{ij}+\zeta^{-ij}, \quad&
\chi_j(b) &= 0.
\end{alignat*}
\end{enumerate}
\end{prop}

Since $a^i$ is a rotation through an angle of $2 \pi/m$ and $a^i b$ is a reflection, 
the values of the class function $\hat{\varphi}$ defined in (\ref{eq:def_phiqu}) are given by
$$
\hat{\varphi}(a^i)
 =
\frac{ ( 1 - \zeta^i u ) ( 1 - \zeta^{-i} u ) }
     { ( 1 - \zeta^i q ) ( 1 - \zeta^{-i} q ) },
\quad
\hat{\varphi}(a^i b)
 =
\frac{ (1-u)(1+u) }
     { (1-q)(1+q) }.
$$
The multiplicities $\hat{m}^\chi$ in $\hat{\varphi}$ can be computed as follows:

\begin{prop}
\label{prop:decomp_D}
\begin{enumerate}
\item[(1)]
If $m$ is even, then we have
\begin{align*}
\hat{m}^{\xi_0}(u) &= \frac{ (1 - u q) (1 - u q^{m-1}) }{ (1-q^2)(1-q^m) },
\\
\hat{m}^{\xi_1}(u) &= q^m \frac{ (1 - u q^{-1}) (1 - u q^{-m+1}) }{ (1-q^2)(1-q^m) },
\\
\hat{m}^{\xi_2}(u) &= q^{m/2} \frac{ (1 - u q) (1 - u q^{-1}) }{ (1-q^2)(1-q^m) },
\\
\hat{m}^{\xi_3}(u) &= q^{m/2} \frac{ (1 - u q) (1 - u q^{-1}) }{ (1-q^2)(1-q^m) },
\\
\hat{m}^{\chi_j}(u) &= (q^j + q^{m-j}) \frac{ (1 - u q) (1 - u q^{-1}) }{ (1-q^2)(1-q^m) }
\quad(1 \le j \le m/2-1).
\end{align*}
\item[(2)]
If $m$ is odd, then we have
\begin{align*}
\hat{m}^{\xi_0}(u) &= \frac{ (1 - u q) (1 - u q^{m-1}) }{ (1-q^2)(1-q^m) },
\\
\hat{m}^{\xi_1}(u) &= q^m \frac{ (1 - u q^{-1}) (1 - u q^{-m+1}) }{ (1-q^2)(1-q^m) },
\\
\hat{m}^{\chi_j}(u) &= (q^j + q^{m-j}) \frac{ (1 - u q) (1 - u q^{-1}) }{ (1-q^2)(1-q^m) }
\quad(1 \le j \le (m-1)/2).
\end{align*}
\end{enumerate}
\end{prop}

\begin{demo}{Proof}
It suffices to show that the right-hand sides of the above formulas satisfy 
the relation $\sum_{\chi \in \Irr(W)} \hat{m}^\chi \chi(w) = \hat{\varphi}(w)$ for $w \in W$.
This can be done by a direct computation, so we leave it to the readers.
\end{demo}

\subsection{%
Proof of Theorem~\ref{thm:main} (A)
}

It is enough to prove (ii) $\implies$ (iii) and (iii) $\implies$ (i).

\begin{demo}{Proof of (ii) $\implies$ (iii)}
We show that, if $\Cat_k(W,q)$ is a polynomial, then $k \equiv \pm1 \bmod m$.
Since the degrees of $D_{2m}$ are $2$ and $m$, we have
$$
\Cat_k(W,q)
 =
\frac{ [k+1]_q [k+m-1]_q }
     { [2]_q [m]_q }.
$$
Suppose that $\Cat_k(W,q)$ is a polynomial in $q$.
Then $\Phi_m(q)$ must appear in the numerator, so we have $m \mid k+1$ or $m \mid k+m-1$, 
which implies $k \equiv -1 \bmod m$ or $k \equiv 1 \bmod m$ respectively.
\end{demo}

\begin{demo}{Proof of (iii) $\implies$ (i)}
We use Lemma~\ref{lem:NtoNq} to prove that, 
if $k \equiv 1$ or $-1 \bmod m$, then $\tilde{m}^\chi_k = \hat{m}^\chi|_{u=q^k}$ 
is a polynomial with nonnegative integer coefficients for all $\chi \in \Irr(W)$.

First consider the multiplicities of $\xi_0$ and $\xi_1$.
By using the formulas given in Proposition~\ref{prop:decomp_D}, we see that
\begin{align*}
\hat{m}^{\xi_0}(qu)
 &=
1
 + \left( q^m + q^2 \frac{ [2m]_q }{ [2]_q } \right) \cdot \frac{ 1-u }{ 1-q^m }
 + q^{2m+2} \frac{ [2m]_q }{ [2]_q } \cdot \frac{ (1-u) (1-uq^{-m}) }{ (1-q^m) (1-q^{2m}) },
\\
\hat{m}^{\xi_0}(q^{m-1} u)
 &=
\frac{ [2m-2]_q }{ [2]_q }
 + q^m \left( \frac{ [2m-2]_q }{ [2]_q } + q^{m-2} \frac{ [2m]_q }{ [2]_q } \right) \cdot \frac{ 1-u }{ 1-q^m }
\\
&\quad\quad
 + q^{4m-2} \frac{ [2m]_q }{ [2]_q } \cdot \frac{ (1-u) (1-uq^{-m}) }{ (1-q^m) (1-q^{2m}) },
\\
\hat{m}^{\xi_1}(qu)
 &=
q^m \cdot \frac{ 1-u }{ 1-q^m}
 + q^{m+2} \frac{ [2m]_q }{ [2]_q } \cdot \frac{ (1-u) (1-uq^{-m}) }{ (1-q^m) (1-q^{2m}) },
\\
\hat{m}^{\xi_1}(q^{m-1}u)
 &=
q^m \frac{ [2m-2]_q }{ [2]_q} \cdot \frac{ 1-u }{ 1-q^m}
 + q^{3m-2} \frac{ [2m]_q }{ [2]_q} \cdot \frac{ (1-u) (1-uq^{-m}) }{ (1-q^m) (1-q^{2m}) }.
\end{align*}
Since $[2m]_q/[2]_q$, $[2m-2]_q/[2]_q \in \Nat[q]$, we conclude that 
$\tilde{m}^{\xi_0}_{pm+1}$, $\tilde{m}^{\xi_0}_{pm+m-1}$, 
$\tilde{m}^{\xi_1}_{pm+1}$ and $\tilde{m}^{\xi_1}_{pm+m-1}$ belong to $\Nat[q]$ 
for all nonnegative integers $p$.

The remaining multiplicities have a factor
$$
M(u) = \frac{ (1-uq)(1-uq^{-1}) }{ (1-q^2) (1-q^m) },
$$
and we see that
\begin{align*}
M(qu)
 &=
\frac{ [m+2]_q }{ [2]_q} \cdot \frac{ 1-u }{ 1-q^m }
 + q^{m-2} \frac{ [2m]_q }{ [2]_q} \cdot \frac{ (1-u) (1-uq^{-m}) }{ (1-q^m) (1-q^{2m}) },
\\
M(q^{m-1}u)
 &=
\frac{ [m-2]_q }{ [2]_q}
 + \left( q^{m-2} \frac{ [2m]_q }{ [2]_q } + q^m \frac{ [m-2]_q }{ [2]_q } \right)
 \cdot \frac{ 1-u }{ 1-q^m }
\\
&\quad\quad
 + q^{3m-2} \frac{ [2m]_q }{ [2]_q } \cdot \frac{ (1-u) (1-uq^{-m}) }{ (1-q^m) (1-q^{2m}) }.
\end{align*}
If $m$ is even, then we have $[m+2]_q/[2]_q$, $[m-2]_q/[2]_q \in \Nat[q]$ 
and $M(q^{pm+1})$, $M(q^{pm+m-1}) \in \Nat[q]$.
Hence $\tilde{m}^\chi_{pm+1}$, $\tilde{m}^\chi_{pm+m-1} \in \Nat[q]$ 
for $\chi = \xi_2$, $\xi_3$ and $\chi_j$.
If $m$ is odd, then we use the relations
\begin{align*}
(q^j + q^{m-j}) \frac{ [m+2]_q }{ [2]_q }
 =
q^j \frac{ [2m-2j+2]_q }{ [2]_q }
+ q^{m-j} \frac{ [2j+2]_q }{ [2]_q },
\\
(q^j + q^{m-j}) \frac{ [m-2]_q }{ [2]_q }
 =
q^j \frac{ [2m-2j-2]_q }{ [2]_q }
+ q^{m-j} \frac{ [2j-2]_q }{ [2]_q }
\end{align*}
to conclude that
$$
\tilde{m}^{\chi_j}_{pm+1} = (q^j+q^{m-j}) M(q^{pm+1}),
\quad
\tilde{m}^{\chi_j}_{pm+m-1} = (q^j+q^{m-j}) M(q^{pm+m-1}).
$$
are polynomials with nonnegative integer coefficients.
\end{demo}

\subsection{%
Proof of Theorem~\ref{thm:main} (C)
}

We prove the implications (iii') $\implies$ (v) and (iv) $\implies$ (iii').

\begin{demo}{Proof of (iii') $\implies$ (v)}
Since $\varphi_1$ is the trivial character, we may assume $k \ge m-1$.

First we consider the case where $m$ is even, and assume that $k^2 \equiv 1 \bmod 2m$.
If we denote by $\eta_1$ and $\eta_2$ the permutation characters on $D_{2m}/\langle b \rangle$ 
and $D_{2m}/\langle ab \rangle$ respectively, then we have
\begin{equation}
\label{eq:decomp_D1}
\varphi_k
 =
\triv + \frac{k-1}{2} \eta_1 + \frac{k-1}{2} \eta_2 + \frac{(k-1)(k-m+1)}{2m} \eta_{\text{reg}},
\end{equation}
where $\eta_{\text{reg}}$ is the character of the regular representation.
Since $k^2 \equiv 1 \bmod 2m$ and $m$ is even, we see that $k$ is odd and
$$
\frac{(k-1)(k-m+1)}{2}
 =
\frac{k^2-1}{2m} - \frac{k-1}{2}
$$
is an integer.
Also, since $k \ge m-1$, the coefficients in (\ref{eq:decomp_D1}) are nonnegative.

Next we consider the case where $m$ is odd, and assume that $k^2 \equiv 1 \bmod m$.
If we denote $\eta_1$ the permutation character on $D_{2m}/\langle b \rangle$, 
then we have
\begin{equation}
\label{eq:decomp_D2}
\varphi_k
 =
\triv + (k-1) \eta_1 + \frac{(k-1)(k-m+1)}{2m} \eta_{\text{reg}}.
\end{equation}
Since $k^2 \equiv 1 \bmod m$ and $m$ is odd, we see that
$(k^2-1)/m$ have the different parity to $k$ and 
$$
\frac{(k-1)(k-m+1)}{2m}
 =
\frac{1}{2}
\left( \frac{k^2-1}{m} - k + 1 \right)
$$
is an integer.
Also, since $k \ge m-1$, the coefficients in (\ref{eq:decomp_D2}) are nonnegative.
\end{demo}

\begin{demo}{Proof of (iv) $\implies$ (iii')}
First we consider the case where $m$ is even.
Then, by using Proposition~\ref{prop:decomp_D}, we have
\begin{align*}
\varphi_k
 &=
\frac{(k+1)(k+m-1)}{2m} \xi_0
+
\frac{(k-1)(k-m+1)}{2m} \xi_1
+
\frac{(k+1)(k-1)}{2m} \xi_2
+
\frac{(k+1)(k-1)}{2m} \xi_3
\\
&\quad+
\sum_{j=1}^{m/2-1} \frac{(k+1)(k-1)}{m} \chi_j.
\end{align*}
Suppose that $\varphi_k$ is the character of some representation of $W$.
Since the multiplicity of $\xi_2$ is an integer, we have $k^2 \equiv 1 \bmod 2m$.
Also the multiplicity of $\xi_1$ is nonnegative, we have $k \le 1$ or $k \ge m-1$.

Similarly we can prove the case where $m$ is odd, by using the decomposition
$$
\varphi_k
 =
\frac{(k+1)(k+m-1)}{2m} \xi_0
+
\frac{(k-1)(k-m+1)}{2m} \xi_1
+
\sum_{j=1}^{(m-1)/2} \frac{(k+1)(k-1)}{m} \chi_j.
$$
\end{demo}

\section{%
Exceptional groups
}

In this section we explain how to prove our main result for the exceptional groups 
with a help of computer.

\subsection{%
Proof of (ii) $\implies$ (iii)
}

We use the criterion for polynomiality in Lemma~\ref{lem:pol_Catalan}.
We can prove the following lemma by using a computer.

\begin{lemma}
\label{lem:qCatalan_exceptional}
Let $W$ be an exceptional complex reflection group with 
degrees $(d_1, \dots, d_r)$ and codegrees $(d^*_1, \dots, d^*_r)$.
For a positive integer $d$ and a positive integer $k$, we put
$$
N_k(d) = \# \{ i : d \mid (k+d_i-1) \},
\quad
N^*_k(d) = \# \{ i : d \mid (k-d_i^*-1) \},
$$
and
$$
D(d) = \# \{ i : d \mid d_i \},
\quad
T = \bigcup_{i=1}^r \{ d : d \mid d_i \}.
$$
Then we have
\begin{enumerate}
\item[(1)]
The following are equivalent for a positive integer $k$:
\begin{enumerate}
\item[(i)]
$N_k(d) \ge D(d)$ for all $d \in T$.
\item[(ii)]
$k$ satisfies the condition in Table~\ref{tab:condition},
except for $W = G_{13}$ and $G_{15}$.
In these exceptions, the condition is given in Table~\ref{tab:qCat}.
\begin{table}[htbp]
\caption{Condition in Lemma~\ref{lem:qCatalan_exceptional} (ii)}
\label{tab:qCat}
$$
\begin{array}{c|c}
\text{group} & \text{condition on $k$} \\
\hline
W = G_{13} & k \equiv 1, 5 \bmod 12 \\
W = G_{15} & k \equiv 1 \bmod 12
\end{array}
$$
\end{table}
\end{enumerate}
\item[(2)]
The following are equivalent for a positive integer $k$:
\begin{enumerate}
\item[(i${}^*$)]
$N^*_k(d) \ge D(d)$ for all $d \in T$.
\item[(ii${}^*$)]
$k$ satisfies the condition in Table~\ref{tab:condition}.
\end{enumerate}
\end{enumerate}
\end{lemma}

\begin{demo}{Proof}
The conditions given in (ii) and (ii${}^*$) are 
of the form ``$k \bmod H \in K$'' for some integer $H$ and 
a subset $K \subset \{ 1, 2, \dots, H \}$, 
where $k \bmod H$ is the remainder of $k$ divided by $H$.
Let $L$ be the least common multiple of $d_1, \dots, d_r$ and $H$.
Then it is easy to see that 
$N_k(d) = N_{L+k}(d)$ and $N^*_k(d) = N^*_{L+k}(d)$.
And $k$ satisfies the condition (ii) (resp. (ii${}^*$)) if and only if 
$k+L$ satisfies (ii) (resp. (ii${}^*$)).
Hence it is enough to show that
\begin{align*}
\{ k \in [1, L] : \text{$N_k(d) \ge D(d)$ for all $d \in T$} \}
 &=
\{ k \in [1, L] : \text{$k$ satisfies (ii)} \},
\\
\{ k \in [1, L] : \text{$N^*_k(d) \ge D(d)$ for all $d \in T$} \}
 &=
\{ k \in [1, L] : \text{$k$ satisfies (ii${}^*$)} \}.
\end{align*}
These equalities can be checked by using a computer.
\end{demo}

Now we are ready to prove (ii) $\implies$ (iii) of Theorem~\ref{thm:main}.

\begin{demo}{Proof of (ii) $\implies$ (iii)}
If follows from Lemmas~\ref{lem:pol_Catalan} and \ref{lem:qCatalan_exceptional} that, 
if $\Cat_k(W,q)$ is a polynomial in $q$, then $k$ satisfies the condition (ii) of 
Lemma~\ref{lem:qCatalan_exceptional}, 
and that,
if $\Cat^*_k(W,q)$ is a polynomial in $q$, then $k$ satisfies the condition (ii${}^*$),
or $k = d^*_i + 1$ for some $i$.

If $k = d^*_i + 1$ does not satisfy the condition (ii${}^*$), then 
the pair $(W,k)$ is one of the following:
$$
(W, k)
 =
(G_{25},4), 
\quad
(G_{33}, 9),
\quad
(G_{33}, 15),
\quad
(G_{35}, 4),
\quad
(G_{35}, 8).
$$
In each of these cases, we can see that $k$ does not satisfy 
the condition (ii) of Lemma~\ref{lem:Catalan_exceptional}.
So we conclude that, if $\Cat_k(W,q)$ and $\Cat^*_k(W,q)$ are both polynomials, 
then $k$ satisfies both (ii) and (ii${}^*$) of Lemma~\ref{lem:qCatalan_exceptional}, 
hence the condition (iii) of Theorem~\ref{thm:main}.
\end{demo}

\subsection{%
Proof of (iii) $\implies$ (i)
}

In order to prove (iii) $\implies$ (i), we need to compute explicitly 
the multiplicities $\hat{m}^\chi$ in $\hat{\varphi}$.
If $\{ w_1, \dots, w_n \}$ is a complete set of representatives of conjugacy classes of $W$, then
the multiplicities are given by
$$
\hat{m}^\chi
 =
\frac{ 1 }{ \# W }
\sum_{w \in W} \hat{\varphi}(w) \overline{\chi(w)}
 =
\sum_{i=1}^n \frac{1}{\# Z(w_i)} \hat{\varphi}(w_i) \overline{\chi(w_i)},
$$
where $Z(w_i)$ is the centralizer of $w_i$ in $W$ and 
$\overline{\chi(w)}$ is the complex conjugate of $\chi(w)$.
By using the program GAP \cite{GAP3} together with CHEVIE \cite{CHEVIE}, 
we can compute 
\begin{itemize}
\item
a complete set of representatives of the conjugacy classes of $W$,
\item
matrix representations of $w \in W$ in the reflection representation,
\item
the character table of $W$,
\item
the size of the centralizer of each representative.
\end{itemize}
Hence we obtain explicit formulas for $\hat{m}^\chi$.
Then we can use Lemma~\ref{lem:NtoNq} to show that $\tilde{m}^\chi_k 
= \hat{m}^\chi|_{u = q^k}$ is a polynomial with nonnegative integer coefficients 
provided $k$ satisfies the condition (iii) (see Table~\ref{tab:condition}) 
of Theorem~\ref{thm:main}.

\begin{remark}
It is known \cite[Corollary~9.39]{LT} 
that the coordinate ring $S = \Comp[V]$ of the reflection representation $V$ 
is decomposed into the tensor product of the invariant subring $S^W$ 
and the space of $W$-harmonic polynomials $\mathcal{H}$.
By using this decomposition, we can show that $\det_V(1 - qw)$ divides 
$\prod_{i=1}^r (1 - q^{d_i})$.
Hence
$$
\prod_{i=1}^r (1 - q^{d_i}) \cdot \hat{m}^\chi
 =
\sum_{i=1}^n
 \frac{1}{\# Z(w_i)} 
 \det{}_V ( 1 - u w )
 \cdot
 \frac{ \prod_{i=1}^r (1 - q^{d_i}) }
      { \det_V ( 1 - q w) }
 \cdot
 \overline{\chi(w_i)}
$$
is a polynomial in $q$ and $u$.
By designing the program in terms of $\prod_{i=1}^r (1 - q^{d_i}) \cdot \hat{m}^\chi$ 
instead of $\hat{m}^\chi$, 
we can avoid rational functions to make the computation faster.
\end{remark}

\subsection{%
Proof of (iii) $\implies$ (v)
}

We use the Orlik--Solomon formula (\ref{eq:OS}) and the data given in \cite{OS1, OS2}.
Then, by using Lemma~\ref{lem:NtoN}, we can prove that 
the coefficients $\chi(\mathcal{A}_j,k)/[N_W(W_j):W_j]$ are nonnegative integers 
provided $k$ satisfies the condition (iii) of Theorem~\ref{thm:main}.

\subsection{%
Proof of (iv) $\implies$ (iii)
}

It follows from Proposition~\ref{prop:Catalan=mult} that 
the Catalan numbers $\Cat_k(W,1)$ and $\Cat^*_k(W,1)$ are the multiplicities of 
$\triv$ and $\det_V$ in $\varphi_k$ respectively.
Hence the proof of (iv) $\implies$ (iii) follows from the following Lemma:

\begin{lemma}
\label{lem:Catalan_exceptional}
Let $W$ be an exceptional complex reflection group.
Then  we have
\begin{enumerate}
\item[(1)]
The following are equivalent for a positive integer $k$:
\begin{enumerate}
\item[(i)]
$\Cat_k(W,1)$ is an integer.
\item[(ii)]
$k$ satisfies the condition given in Table~\ref{tab:condition}, 
except for $W = G_{13}$, $G_{15}$, $G_{25}$, $G_{33}$, $G_{35}$ and $G_{36}$.
In these exceptions, the condition is given in Table~\ref{tab:Cat}.
\begin{table}[htbp]
\caption{Condition in Lemma~\ref{lem:Catalan_exceptional} (ii)}
\label{tab:Cat}
$$
\begin{array}{c|c}
\text{group} & \text{condition} \\
\hline
G_{13} & k \equiv 1, 5 \bmod 12 \\
G_{15} & k \equiv 1 \bmod 12 \\
G_{25} & k \equiv 1 \bmod 6 \quad\text{or}\quad k \equiv 16 \bmod 24 \\
G_{33} & k \equiv 1 \bmod 6 \quad\text{or}\quad k \equiv 45, 51 \bmod 54 \\
G_{35} & k \equiv 1, 5 \bmod 6 \quad\text{or}\quad k \equiv 28,56,88,92 \bmod 96 \\
G_{36} & k \equiv 1, 5 \bmod 6 \quad\text{or}\quad k \equiv 153 \bmod 162
\end{array}
$$
\end{table}
\end{enumerate}
\item[(2)]
The following are equivalent for a positive integer $k$:
\begin{enumerate}
\item[(i${}^*$)]
$\Cat^*_k(W,1)$ is an integer.
\item[(ii${}^*$)]
$k$ satisfies the condition given in Table~\ref{tab:condition}, 
except for $W = G_{25}$, $G_{33}$, $G_{35}$, and $G_{36}$.
In these exceptions, the condition is given in Table~\ref{tab:Cat*}.
\begin{table}[htbp]
\caption{Condition in Lemma~\ref{lem:Catalan_exceptional} (ii${}^*$)}
\label{tab:Cat*}
$$
\begin{array}{c|c}
\text{group} & \text{condition} \\
\hline
G_{25} & k \equiv 1 \bmod 6 \quad\text{or}\quad k \equiv 4 \bmod 24 \\
G_{33} & k \equiv 1 \bmod 6 \quad\text{or}\quad k \equiv 9, 15 \bmod 54 \\
G_{35} & k \equiv 1, 5 \bmod 6 \quad\text{or}\quad k \equiv 4,8,40,68 \bmod 96 \\
G_{36} & k \equiv 1, 5 \bmod 6 \quad\text{or}\quad k \equiv 9 \bmod 162
\end{array}
$$
\end{table}
\end{enumerate}
\end{enumerate}
\end{lemma}

\begin{demo}{Proof}
We note that
$$
c(k) = \Cat_k(W,1) = \prod_{i=1}^r \frac{k+d_i-1}{d_i},
\quad
c^*(k) = \Cat^*_k(W,1) = \prod_{i=1}^r \frac{k-d~*_i-1}{d_i}
$$
are polynomials in $k$.
And the conditions in (ii) and (ii${}^*$) are of the form 
``$k \bmod H \in K$'' for some positive integer $H$ and a subset $K \subset \{ 1, \dots, H \}$.
Then we can find a positive integer $L$ satisfying
\begin{enumerate}
\item[(a)]
$L$ is a multiple of $H$, and
\item[(b)]
$c(t+L) - c(t)$ and $c^*(t+L) - c^*(t)$ map nonnegative integers to nonnegative integers.
\end{enumerate}
For example, we can take $L = \# W$.
(We can use Lemma~\ref{lem:NtoN} to check the condition (b).)
Then it is enough to show that
\begin{align*}
\{ k \in [1, L] : c(k) = \Cat_k(W,1) \in \Nat \}
 &=
\{ k \in [1, L] : \text{$k$ satisfies (ii)} \},
\\
\{ k \in [1, L] : c^*(k) = \Cat^*_k(W,1) \in \Nat \}
 &=
\{ k \in [1, L] : \text{$k$ satisfies (ii${}^*$)} \}.
\end{align*}
These equalities can be checked by using a computer.
\end{demo}

This complete the proof of Theorem~\ref{thm:main} for all irreducible finite complex reflection groups.
\section{%
Concluding remarks
}

\subsection{%
Uniform proof
}

We have proved Theorem~\ref{thm:main} in a case-by-case manner.
The conditions (i) and (ii) of Theorem~\ref{thm:main} are stated in a uniform way 
(i.e., do not use the classification of irreducible complex reflection groups).
So it is natural to seek a classification-free proof.

\begin{problem}
Find a uniform proof to the equivalence between (i) and (ii) of Theorem~\ref{thm:main}.
\end{problem}

If $W$ is a Weyl group, then the condition (iii) is equivalent to saying that 
$k$ is ``very good'' in the sense of \cite{Sommers}. 
A positive integer $k$ is called ``very good'' for $W$ 
if $k$ is relatively prime to the coefficients of simple roots in the highest root 
and to the index $[Q^\vee : Q]$, 
where $Q$ and $Q^\vee$ is the root lattice and the coroot lattice respectively.

\begin{problem}
Find a classification-free description of the condition (iii) in Theorem~\ref{thm:main}.
\end{problem}

\subsection{%
Combinatorial/algebraic models
}

In Theorem~\ref{thm:main}, we have proved that, if $k$ satisfies the condition (iii) 
(or (iii') for dihedral groups), 
then $\tilde{\varphi}_k$ is the graded character of a representation 
and $\varphi_k$ is a permutation character.
However, our proof does not give any constructions of such representations.
We know \cite{Sommers} that, if $W$ is a Weyl group, then $Q/kQ$ affords the character $\varphi_k$ 
provided $k$ is very good.
Also, for a Coxeter group $W$ and $k \equiv 1 \bmod h$ (Fuss cases), 
\cite{ARR} and \cite{Rhoades} give a combinatorial model for $\varphi_k$.
So it is a natural problem to generalize their constructions to other complex reflection groups 
and to general $k$ satisfying the condition (iii) in Theorem~\ref{thm:main}.

\begin{problem}
\begin{enumerate}
\item[(1)]
Give a combinatorial construction of a $W$-set which gives the permutation character $\varphi_k$.
\item[(2)]
Give an algebraic/combinatorial construction of a graded $W$-module affording the graded character 
$\tilde{\varphi}_k$.
\end{enumerate}
\end{problem}

\subsection{%
Unimodality
}

In Proposition~\ref{prop:skd_Nq}, we have proved that
$$
f_\lambda(q)
 = 
\frac{ s_\lambda(1, q, \dots, q^{k-1}) }
     { [k]_q/[d]_q }
$$
is a polynomial in $q$ with nonnegative integer coefficients, 
where $\lambda$ is a partition of $n$ and $d = \gcd(n,k)$.
Recall that a finite sequence $(a_0, a_1, \dots, a_m)$ is unimodal if there is an index $p$ such that
$$
a_0 \le a_1 \le \dots \le a_{p-1} \le a_p \ge a_{p+1} \ge \dots \ge a_m.
$$
If $n$ is a multiple of $k$, i.e., $\gcd(n,k) = k$, then it is well-known 
(see \cite[I.8 Example~4]{Macdonald} for example) that 
the coefficients of $f_\lambda(q) = s_\lambda(1, q, \dots, q^{k-1})$ 
form a unimodal sequence.
Also, if $\gcd(n,k) = 1$, then a certain finite dimensional module over a rational Cherednik algebra 
associated to $\Sym_n$ provides a $W \times {\mathfrak{sl}}_2$ module with character $\tilde{\varphi}_k$.
(See \cite{BEG} for example.)
Hence the polynomials $f_\lambda(q)$ are ${\mathfrak{sl}}_2$ characters.
These results together with a computer experiment suggest the following conjecture.

\begin{conjecture}
Let $k$ and $n$ be positive integer and put $d = \gcd(n,k)$.
Given a partition $\lambda$ of $n$ with $l(\lambda) \le k$, we write
$$
\frac{ s_\lambda(1, q, \dots, q^{k-1}) }
     { [k]_q/[d]_q }
=
\sum_{i \ge 0} a_i q^i.
$$
Then the sequences $(a_0, a_2, a_4, \dots)$ and $(a_1, a_3, a_5, \dots)$ are both unimodal.
\end{conjecture}

\begin{remark}
The whole sequence $(a_0, a_1, a_2, a_3, \dots)$ is not unimodal in general.
For example, if $n=2$, $k=3$ and $\lambda = (2)$, then
$$
\frac{ s_{(2)}(1, q, q^2) }
     { [3]_q }
 =
1 + q^2
$$
does not have unimodal coefficients.
\end{remark}
\appendix
\section{%
Divisibility of the Stirling numbers
}

In this appendix, we give a proof to Lemma~\ref{lem:Stirling} which asserts that 
the signless Stirling number $c(n,j)$ of the first kind is
 divisible by $\binom{n}{2}$ if $n-j$ is odd.
Recall that $c(n,j)$ is equal to 
the number of permutations of $n$ letters with exactly $j$ cycles.
If $w \in \Sym_n$ has the cycle type $\lambda$, then we have 
$n - l(\lambda) = \sum_{i \ge 1} (i-1) m_i(\lambda)$, 
where $m_i(\lambda)$ is the multiplicity of $i$ in $\lambda$.
Hence, if $n - l(\lambda)$ is odd, then $\lambda$ has an even part with odd multiplicity.
Now Lemma~\ref{lem:Stirling} follows from the following proposition.
(According to \cite[Exercise~91 (c)]{EC2supp}, J.~Bruns proved the following result, 
but it seems that he never published his proof.)

\begin{prop}
\label{prop:refinedStirling}
Let $\lambda$ be a partition of $n$ and $C_\lambda$ the corresponding conjugacy class of $\Sym_n$.
If $\lambda$ has an even part with odd multiplicity, 
then $\# C_\lambda$ is divisible by $\binom{n}{2}$.
\end{prop}

\begin{demo}{Proof}
Suppose that $\lambda$ has an even part $2f$ with odd multiplicity $l = m_{2f}(\lambda)$.
The size of the conjugacy class $C_\lambda$ is given by
$\# C_\lambda = n! / z_\lambda$, where $z_\lambda = \prod_{i \ge 1} i^{m_i(\lambda)} m_i(\lambda)!$.
Let $\mu$ be the partition obtained by removing $l$ parts equal to $2f$ from $\lambda$.
Then we have
$$
|\mu| = n - 2fl,
\quad
z_\lambda = z_\mu \cdot (2f)^l \cdot l!.
$$
Hence we have
$$
\frac{ 1 }{ \binom{n}{2} }
\# C_\lambda
 =
\binom{ n-2 }{ 2fl -2 }
\cdot
\# C_\mu
\cdot
\frac{ (2fl-2) ! }
     { f^l \cdot (2l-2)! }
\cdot
\frac{ (2l-2)! }
     { 2^{l-1} \cdot l! }.
$$
If $f=1$, then $(2fl-2) !/ f^l \cdot (2l-2)! = 1$.
If $f \ge 2$, then $2l-2 \le fl-2$ and the interval $[fl-1,2fl-2]$ contains 
$l$ multiples of $f$, 
so we see that $(2fl-2)!/(2l-2)!$ is divisible by $f^l$.
It remains to show that the last factor $(2l-2)! / 2^{l-1} \cdot l!$ is an integer.

Let $\nu_2(x)$ denote the $2$-adic valuation of a rational number $x$.
Then we have
$$
\nu_2 \left( \frac{ (2l-2)! }{ l! } \right)
 =
\sum_{i \ge 1} \left\lfloor \frac{ 2l-2 }{ 2^i } \right\rfloor
-
\sum_{i \ge 1} \left\lfloor \frac{ l }{ 2^i } \right\rfloor,
$$
where $\lfloor x \rfloor$ denotes the largest integer not exceeding $x$.
Since $l$ is an odd integer, we have $\lfloor (2l-2) / 2^{i+1} \rfloor
 = \lfloor (l-1) / 2^i \rfloor = \lfloor l/2^i \rfloor$.
Hence we have
$$
\nu_2 \left( \frac{ (2l-2)! }{ l! } \right)
 =
\left\lfloor \frac{ 2l-2 }{ 2 } \right\rfloor
 =
l-1,
$$
and $(2l-2)!/l!$ is exactly divisible by $2^{l-1}$.
\end{demo}
\section*{%
Acknowledgments
}

The authors would like to thank Professors C.~Krattenthaler 
and R.~Stanley for their helpful comments.
The second author is partially supported by the JSPS Grants-in-Aid 
for Scientific Research No.~24340003.


\end{document}